\documentclass[twoside]{article} 
\usepackage[accepted]{aistats2018}
\usepackage[round]{natbib}

\usepackage[utf8]{inputenc} %
\usepackage[T1]{fontenc}    %
\usepackage{hyperref}       %
\usepackage{url}            %
\usepackage{booktabs}       %
\usepackage{amsfonts}       %
\usepackage{amsmath}
\usepackage{nicefrac}       %
\usepackage{microtype}      %
\usepackage{algcompatible} %
\usepackage{algorithm}
\usepackage{caption}
\usepackage{subcaption}
\usepackage{color}
\usepackage{nccmath}
\definecolor{mydarkblue}{rgb}{0,0.08,0.45}

\hypersetup{ %
    pdftitle={Frank-Wolfe Splitting via Augmented Lagrangian Method},
    pdfkeywords={},
    pdfborder=0 0 0,
  pdfpagemode=UseNone,
    colorlinks=true,
    linkcolor=mydarkblue, %
    citecolor=mydarkblue, %
    filecolor=mydarkblue, %
    urlcolor=mydarkblue, %
    pdfview=FitH}

\usepackage{bm}
\usepackage{mathdots}
\usepackage{amsmath}
\usepackage{amsthm}
\usepackage{amssymb}
\usepackage{sectsty}
\usepackage{graphicx}
\usepackage{color}
\usepackage{todonotes}
\usepackage{mdframed}
\usepackage{stmaryrd}
\usepackage{hyperref}
\usepackage{pgfkeys}
\usepackage{pgfplots}
\usetikzlibrary{shapes,arrows,calc}

\usetikzlibrary{patterns}
\usetikzlibrary{positioning}
\usepackage{empheq}
\definecolor{myblue}{HTML}{D2E4FC}
\newcommand*\mybluebox[1]{%
\colorbox{myblue}{\hspace{1em}#1\hspace{1em}}}

\usepackage{wrapfig}

\newcommand{\N}{\mathbb{N}}
\newcommand{\R}{\mathbb{R}}
\newcommand{\LL}{\mathcal{L}}
\newcommand{\X}{\mathcal{X}}
\newcommand{\Y}{\mathcal{Y}}
\newcommand{\x}{\bm{x}}
\newcommand{\y}{\bm{y}}
\newcommand{\innerProd}[2]{\left\langle #1 , #2 \right\rangle}
\newcommand{\innerProdCompressed}[2]{\langle #1 , #2 \rangle}
\DeclareMathOperator*{\argmin}{\arg\min}
\DeclareMathOperator*{\argmax}{\arg\max}
\newmdtheoremenv{alg}{Algorithm}
\newmdtheoremenv{theo}{Theorem}

\newtheorem{theorem}{Theorem}
\newtheorem{lemma}{Lemma}
\newtheorem{proposition}{Proposition}

\newtheorem{corollary}{Corollary}

\newtheorem{assumption}{Assumption}

\DeclareMathOperator*{\minimize}{minimize}

\newtheorem*{rep@theorem}{\rep@title}
\newcommand{\newreptheorem}[2]{%
\newenvironment{rep#1}[1]{%
 \def\rep@title{#2 \ref{##1}}%
 \begin{rep@theorem}}%
 {\end{rep@theorem}}}

\newreptheorem{lemma}{Lemma'}
\newreptheorem{proposition}{Proposition'}
\newreptheorem{theorem}{Theorem'}
\newreptheorem{assumption}{Assumption'}

\DeclareMathOperator{\diam}{diam}

\DeclareMathOperator*{\diag}{diag}
\DeclareMathOperator*{\dist}{dist}

\DeclareMathOperator*{\relint}{\mathop{Relint}}
\DeclareMathOperator*{\cone}{\mathop{Cone}}
\DeclareMathOperator*{\Span}{\mathop{Span}}

\newcommand{\bigO}{O}

\newcommand{\stepsize}{\gamma}
\newcommand{\stepmax}{\stepsize_{\textnormal{\scriptsize max}}} %

\newcommand{\s}{\bm{s}}
\newcommand{\dd}{\bm{d}}
\newcommand{\uu}{\bm{u}}
\newcommand{\vv}{\bm{v}} %

\DeclareMathOperator*{\lmo}{\textsc{Lmo}}

\newcommand{\Vertices}{\mathcal{A}} %
\newcommand{\Coreset}{\mathcal{S}}

\renewcommand{\r}{\bm{r}}

\newcommand{\C}{\mathcal{C}}

\newcommand{\ADMM}{\textsc{ADMM}}
\newcommand{\AL}{\textsc{ALM}}
\newcommand{\LMO}{\textsc{LMO}}
\newcommand{\GDMM}{\textsc{GDMM}}
\newcommand{\FWAL}{\textsc{FW-AL}}
\newcommand{\FW}{\textsc{FW}}

\begin{document}

\twocolumn[

\aistatstitle{Frank-Wolfe Splitting via Augmented Lagrangian Method}

\aistatsauthor{ Gauthier Gidel \And Fabian Pedregosa \And  Simon Lacoste-Julien }
\aistatsaddress{
 MILA, DIRO Université de Montréal
 \And  UC Berkeley \& ETH Zurich
 \And MILA, DIRO Université de Montréal

 }
 ]

\begin{abstract}
  Minimizing a function over an intersection of convex sets is an important task in optimization that is often much more challenging than minimizing it over each individual constraint set.
  While traditional methods such as Frank-Wolfe (\FW) or proximal gradient descent assume access to a linear or quadratic oracle on the intersection,
  splitting techniques take advantage of the structure of each sets, and only require access to the oracle on the individual constraints. In this work, we develop and analyze the Frank-Wolfe Augmented Lagrangian (\FWAL) algorithm, a method for minimizing a smooth function over convex compact sets related by a ``linear consistency'' constraint that only requires access to a linear minimization oracle over the individual constraints. It is based on the Augmented Lagrangian Method (\AL), also known as Method of Multipliers,
  but unlike most existing splitting methods, it only requires access to linear (instead of quadratic) minimization oracles. We use recent advances in the analysis of Frank-Wolfe and the alternating direction method of multipliers algorithms to prove a sublinear convergence rate for \FWAL\ over general convex compact sets and a linear convergence rate over polytopes.
\end{abstract}

\section{Introduction} %
\label{sec:introduction}

The Frank-Wolfe (\FW) or conditional gradient algorithm has seen an impressive revival in recent years, notably due to its very favorable properties for the optimization of sparse problems~\citep{jaggi_revisiting_2013} or over structured constraint sets~\citep{lacoste-julien_global_2015}. 
This algorithm assumes knowledge of a linear minimization oracle (\LMO) over the set of constraints. 
This oracle is inexpensive to compute for sets such as the $\ell_1$ or trace norm ball. However, inducing complex priors often requires to consider \emph{multiple} constraints, leading to a constraint set formed by the intersection of the original constraints. Unfortunately, evaluating the \LMO\ over this intersection may be challenging even if the \LMO s\ on the individual sets are inexpensive.

The problem of minimizing over an intersection of convex sets is pervasive in machine learning and signal processing. For example, one can seek for a matrix that is both sparse and low rank by constraining the solution to have \emph{both} small $\ell_1$ and trace norm~\citep{richard_estimation_2012}
or find a set of brain maps which are both sparse and piecewise constant by constraining both the $\ell_1$ and total variation pseudonorm~\citep{gramfort2013identifying}. Furthermore, some challenging optimization problems such as multiple sequence alignment are naturally expressed over an intersection of sets~\citep{yen_convex_2016} or more generally as a linear relationship between these sets~\citep{huang_greedy_2017}.

The goal of this paper is to describe and analyze \FWAL, an optimization method that solves convex optimization problems subject to multiple constraint sets, assuming we have access to a \LMO\ on each of the set.

\paragraph{Previous work.} %
\label{par:related_work_}
The vast majority of methods proposed to solve optimization problems over an intersection of sets rely on the availability of a projection operator onto each set (see e.g. the recent reviews \citep{glowinski2017splitting,ryu_primer_2016}, which cover the more general proximal splitting framework).
One of the most popular algorithms in this framework is the alternating direction method of multipliers (\ADMM), proposed by \citet{glowinski_sur_1975}, studied by \citet{gabay_dual_1976}, and revisited many times; see for instance \citep{boyd_distributed_2011,yan_self_2016}.
On some cases, such as constraints on the trace norm~\citep{cai2010singular} or the latent group lasso~\citep{obozinski2011group}, the projection step can be a time-consuming operation, while the Frank-Wolfe \LMO\ is much cheaper in both cases. Moreover, for some highly structured polytopes such as those appearing in alignment constraints~\citep{alayrac16unsupervised} or Structured SVM~\citep{lacostejulien:hal-00720158}, there exists a fast and elegant dynamic programming algorithm to compute the \LMO, while there is no known practical algorithm to compute the projection.
Hence, the development of splitting methods that use the \LMO\ instead of the proximal operator is of key practical interest. 
\citet{yurtsever2015universal} proposed a general algorithm (UniPDGrad) based on the Lagrangian method which, with some work, can be reduced to a splitting method using only LMO as a particular case. We develop the comparison with FW-AL in App.~\ref{sub:comparison_with_unipdgrad}.

Recently, \citet{yen_convex_2016} proposed a \FW\ variant for objectives with a linear loss function over an intersection of polytopes named Greedy Direction Method of Multipliers (\GDMM). A similar version of \GDMM\ is also used in \citep{yen_dual_2016,huang_greedy_2017} to optimize a function over a Cartesian product of spaces related to each other by a linear consistency constraint.
The constraints are incorporated through the augmented Lagrangian method and its convergence analysis crucially uses recent progress in the analysis of \ADMM\ by~\citet{hong_linear_2017}. Nevertheless, we argue in Sec.~\ref{sub:properties_of_the_dual_function_} that there are technical issues in these analysis since some of the properties used have only been proven for \ADMM\ and do not hold in the context of \GDMM.
Furthermore, even though \GDMM\ provides good experimental results in these papers, the practical applicability of the method to other problems is dampened by overly restrictive assumptions: the loss function is required to be linear or quadratic, leaving outside loss functions such as logistic regression, and the constraint needs to be a polytope, leaving outside domains such as the trace norm ball.

\paragraph{Contributions.} %
\label{par:contributions_}
Our main contribution is the development of a novel variant of \FW\ for the optimization of a function over
product of spaces related to each other by a linear consistency constraint and its rigorous analysis. We name this method Frank-Wolfe via Augmented Lagrangian method (\FWAL).
With respect to \citet{yen_convex_2016,yen_dual_2016,huang_greedy_2017},
our framework generalizes \GDMM\ by providing a method to optimize a general class of functions over an intersection of an arbitrary number of compact sets, which are \emph{not} restricted to be polytopes.
Moreover, we argue that the previous proofs of convergence were incomplete: in this paper, we prove a new challenging technical lemma providing a growth condition on the augmented dual function which allows us to fix the missing parts.

We show that \FWAL\ converges for any smooth objective function. We prove that a standard gap measure converges linearly (i.e., with a geometric rate) when the constraint sets are polytopes, and sublinearly for general compact convex sets. We also show that when the function is strongly convex, the sum of this gap measure and the feasibility gives a bound on the distance to the set of optimal solutions.
This is of key practical importance since the applications that we consider (e.g., minimization with trace norm constraints) verify these assumptions.

The paper is organized as follows. In Sec.~\ref{sec:problem_setting}, we introduce the general setting, provide some motivating applications and present the augmented Lagrangian formulation of our problem. In Sec.~\ref{sec:the_algorithm}, we describe the algorithm \FWAL\ and provide its analysis in Sec.~\ref{sec:analysis_of_fw}. Finally, we present illustrative experiments in Sec.~\ref{sec:experiments}.

\section{Problem Setting} %
\label{sec:problem_setting}

We will consider the following minimization problem,
\begin{empheq}[box=\mybluebox]{equation}\tag{OPT}\label{eq:opt}\centering
  \begin{aligned}
  & \quad \quad\minimize_{\x^{(1)}, \ldots, \x^{(k)}}^{\vphantom{i}} \, f(\x^{(1)},\ldots,\x^{(k)}) \;, \\
  & \! \!\!  \text{s.t.} \; \x^{(k)} \in \X_k, \; k \in [K], \; \sum_{k=1}^K A_k \x^{(k)} \!= 0\,, \!\!
  \end{aligned}
\end{empheq}
where $f : \R^m \to \R$ is a convex differentiable function and for $k \in [K]$, $\X_k \subset \R^{d_k}$ are convex compact sets and $A_k$ are matrices of size $d\times d_k$. We will call the constraint $ \sum_{k=1}^K A_k \x^{(k)} = 0$ the \emph{linear consistency constraint}, motivated by the marginalization \emph{consistency} constraints appearing in some of the applications of our framework as described in Sec.~\ref{sec:applications}. One important potential application is the \textbf{intersection of multiple sets}. The simple $K=2$ example can be expressed with $A_1 = I$ and $A_2 = -I$. 
We assume that we have access to the linear minimization oracle $\text{\LMO}_{k}(\r) \in \argmin_{\s \in \X_k} \innerProd{\s}{\r}, \; k \in [K]$. 
We denote by $\X^*$ the set of optimal points of the optimization problem~\eqref{eq:opt} and we assume that this problem is feasible, i.e., the set of solutions is non empty.

\subsection{Motivating Applications} \label{sec:applications}
We now present some motivating applications of our problem setting, including examples where special case versions of \FWAL\ were used.
This previous work provides additional evidence for the practical significance of the \FWAL\ algorithm.
Multiple sequence alignment and motif discovery~\citep{yen_convex_2016}
are problems in which the domain is described as an intersection of alignment constraints and consensus constraints, two highly structured polytopes. The LMO on both sets can be solved by dynamic programing whereas there is no known practical algorithm to project onto.
A factorwise approach to the dual of the structured SVM objective~\citep{yen_dual_2016} can also be cast as constrained problem over a Cartesian product of polytopes related to each other by a linear consistency constraint. As often in structured prediction, the output domain grows exponentially, leading to very high dimensional polytopes. Once again, dynamic programming can be used to compute the linear oracle in structured SVMs at a lower computational cost than the potentially intractable projection.
The algorithms proposed by \citet{yen_convex_2016} and \citet{yen_dual_2016} are in fact a particular instance of \FWAL, where the objective function is respectively linear and quadratic.

Finally, simultaneously sparse ($\ell_1$ norm constraint) and low rank (trace norm constraint) matrices~\citep{richard_estimation_2012} is another class of problems where the constraints consists of an intersection of sets with simple LMO but expensive projection. This example is a novel application of \FWAL{} and is developed in Sec.~\ref{sec:experiments}.

\subsection{Augmented Lagrangian Reformulation} %
\label{sub:augmented_lagrangian_reformulation}
It is possible to reformulate~\eqref{eq:opt} into the problem of finding a saddle point of an augmented Lagrangian~\citep{bertsekas-AL}, in order to split the constraints in a way in which the linear oracle is computed over a product space. We first rewrite~\eqref{eq:opt} as follows:
\begin{equation}\label{eq:opt_reformulation}
\min_{\x^{(k)} \in \X_k, \,k \in [K]} f(\x)
  \quad \text{s.t.} \quad
  M \x = 0 \; ,
\end{equation}
where $\x := \left( \x^{(1)},\ldots,\x^{(K)}\right)$ and $M := \left[ A_1,\ldots,A_k \right]$ is such that,
\begin{equation}\label{eq:form_M}
M\x = 0 \Leftrightarrow \sum_{k=1}^K A_k\x^{(k)} = 0\,.  
\end{equation}
We can now consider the augmented Lagrangian formulation of \eqref{eq:opt_reformulation}, where $\y$ is the dual variable:
\begin{equation}
\label{eq:opt_bis} \tag{OPT2}
\begin{aligned}
   &\minimize_{(\x^{(1)},\ldots,\x^{(K)})} \, \max_{\y \in \R^d} \LL(\x^{(1)},\ldots,\x^{(K)},\y) \\
  & \;\;\text{s.t.} \quad \x^{(k)} \in \X_k, \;\; k\in \{1,\ldots,K\} \; \\
  & \LL(\x,\y) :=  f(\x) + \innerProd{\y}{M\x} + \tfrac{\lambda}{2} \|M\x\|^2 .
\end{aligned}
\end{equation}
We note $\X := {\X_1 \times \cdots \times \X_K} \subset \R^{d_1 + \ldots +d_K} = \R^m$ for notational simplicity. This formulation is the one onto which our algorithm \FWAL{} is applied.

\paragraph{Notation and assumption.} %
\label{par:properties_of_}
In this paper, we denote by $\|\cdot\|$ the $\ell_2$ norm for vectors (resp. spectral norm for matrices) and $\dist(\x, \mathcal C) := \inf_{\x' \in \C} \|\x - \x'\|$ its associated distance to a set. 
We assume that $f$ is $L$-smooth on $\R^m$, i.e., differentiable with $L$-Lipschitz continuous gradient:
\begin{equation}
  \| \nabla f(\x) - \nabla f(\x') \| \leq L\|\x - \x'\| \quad \forall \x, \x' \in  \R^m\,.
\end{equation}
This assumption is standard in convex optimization~\citep{nesterov_introductory_2004}. Notice that the \FW\ algorithm does not converge if the objective function is not at least continuously differentiable~\citep[Example 1]{nesterov_complexity_2016}.
In our analysis, we will also use the observation that $\frac{\lambda}{2}\|M\cdot\|^2$ is generalized strongly convex.\footnote{This notion has been studied by \citet{wang_iteration_2014} and in the Frank-Wolfe framework by \citet{beck_linearly_2016} and \citet{lacoste-julien_global_2015}.} We say that a function $h$ is \emph{generalized strongly convex} when it takes the following general form:
\begin{equation}\label{eq:generalized_strong_conv}
h(\x) := g(A \x) + \innerProd{\bm{b}}{\x}, \quad \forall \x \in \R^m \; ,
\end{equation}
where $A \in \R^{d \times m}$ and $g$ is $\mu_g$-\emph{strongly convex}  w.r.t.\ the Euclidean norm on $\R^d$ with $\mu_g > 0$.
Recall that a $\mu_g$-strongly (differentiable) convex function $g :\R^{d} \to \R$ is one such that, $\forall \x, \x' \in \R^d,$
\begin{equation}\notag
   g(\x) \geq g(\x') + \innerProd{\x - \x'}{\nabla g(\x')} + \frac{\mu_g}{2} \|\x - \x'\|^2 \;.
\end{equation}

\section{\FWAL\ Algorithm} %
\label{sec:the_algorithm}
Our algorithm takes inspiration from both Frank-Wolfe and the augmented Lagrangian method. The augmented Lagrangian method alternates a primal update on $\x$ (approximately) minimizing\footnote{An example of approximate minimization is taking one proximal gradient step, as used for example, in the Linearized ADMM algorithm~\citep{Goldfarb2013, yang2013linearized}. }
the augmented Lagrangian $\LL(\cdot, \y_t)$, with a dual update on $\y$ by taking a gradient ascent step on $\LL(\x_{t+1}, \cdot)$.
The \FWAL\ algorithm follows the general iteration of the augmented Lagrangian method, but with the crucial difference that Lagrangian minimization is replaced by one Frank-Wolfe step on $\LL(\cdot, \y_t)$. The algorithm is thus composed by two loops: an outer loop presented in~\eqref{eq:fwal} and an inner loop noted $\mathcal{FW}$ which can be chosen to be one of the \FW\ step variants described in Alg.~\ref{alg:AFW} or~\ref{alg:FW}.

\paragraph{FW steps.} %
\label{par:the_algorithm}
In \FWAL\, we need to ensure that the $\mathcal{FW}$ inner loop makes sufficient progress. For general sets, we can use one iteration of the classical Frank-Wolfe algorithm with line-search~\citep{jaggi_revisiting_2013} as given in Algorithm~\ref{alg:FW}. When working over polytopes, we can get faster (linear) convergence by taking one \emph{non-drop} step (defined below) of the away-step variant of the FW algorithm (AFW)~\citep{lacoste-julien_global_2015}, as described in Algorithm~\ref{alg:AFW}. Other possible variants are discussed in Appendix~\ref{app:frank_wolfe_algorithms}.
We denote by $\x_{t}$ and $\y_{t}$ the iterates computed by \FWAL\ after $t$ steps and by $\Vertices_{t}$ the set of atoms previously given by the FW oracle (including the initialization point).
\begin{algorithm}[t]
  \caption{Away-step Frank-Wolfe (one non-drop step) : \citep{lacoste-julien_global_2015}}
  \label{alg:AFW}
  \begin{algorithmic}[1]
  \STATE \textbf{input:} ($\x,\,\Coreset,\, \Vertices,\,\varphi$)
  \hfill \emph{($\varphi$ is the objective)}
  \STATE \texttt{drop\_step} $\leftarrow$ \texttt{true} (initialization of the boolean)
  \WHILE{\texttt{drop\_step} = \texttt{true}}
  \STATE $\bm{s} \leftarrow \lmo \left(\nabla \varphi(\x)\right)$
    \STATE $\vv \in \argmax_{\vv \in \Coreset } \left\langle \nabla \varphi(\x), \vv \right\rangle$
    \STATE $g^{FW} \leftarrow \innerProd{\nabla \varphi(\x)}{\x- \bm{s}}$  \hfill \emph{(Frank-Wolfe gap)}
    \STATE $g^{A}  \leftarrow \innerProd{\nabla \varphi(\x)}{\vv - \x}$ \hfill \emph{(Away gap)}
      \IF{ $g^{FW} \geq g^{A}$} \hfill \emph{(FW direction is better)}
      \STATE $\dd \leftarrow  \s - \x$ and $\stepmax \leftarrow 1$
      \ELSE  \hfill \emph{(Away direction is better)}
      \STATE $\dd \leftarrow  \x - \vv$ and $\stepmax \leftarrow \alpha_{\vv} / (1- \alpha_{\vv})$\label{algLine:drop_step}

      \ENDIF
      \STATE Compute $\stepsize \!\in\! \argmin_{\stepsize \in [0,\stepmax]} \textstyle \varphi\left(\x + \stepsize \dd\right)$\label{line:line_search}
      \IF{$\stepsize < \stepmax$} \hfill \emph{(first non-drop step)}
        \STATE \texttt{drop\_step} $\leftarrow$ \texttt{false}
      \ENDIF
      \STATE Update $\x \leftarrow \x + \stepsize \dd$
      \STATE Update $\alpha_{\vv}$ according to \eqref{eq:active_set_decomposition}
      \STATE Update $\Coreset \leftarrow \{\vv \in \Vertices \,\: \mathrm{ s.t. } \,\: \alpha_{\vv} > 0\}$ \hfill \emph{(active set)}
    \ENDWHILE
    \STATE \textbf{return}: $(\x ,\Coreset)$
  \end{algorithmic}
\end{algorithm}
\begin{algorithm}[t]
    \caption{\FW (one step) :~\citep{frank_algorithm_1956}}\label{alg:FW}
    \begin{algorithmic}[1]
      \STATE \textbf{input:} $(\x,\varphi)$
      \STATE Compute $\bm{s} \leftarrow \displaystyle \argmin_{\bm{s} \in \X}\innerProd{\bm{s}}{\nabla \varphi(\x)}$
      \STATE $\gamma \in\argmin_{\gamma \in [0,1] } \varphi(\x + \gamma(\s-\x))$
      \STATE Update $\x \leftarrow (1- \gamma) \x + \gamma \s$
      \STATE \textbf{return:} $\x$
    \end{algorithmic}
  \end{algorithm}
If the constraint set is the convex hull of a set of atoms $\Vertices$, the iterate $\x_{t}$ has a sparse representation as a convex combination of the first iterate and the atoms previously given by the \FW\ oracle. The set of atoms which appear in this expansion with non-zero weight is called the \emph{active set}~$\Coreset_t$. Similarly, since $\y_{t}$ is by construction in the cone generated by $\{M\x_s\}_{s\leq t}$, the iterate $\y_{t}$ is in the span of $M\Vertices_t$, that is, they both have the sparse expansion:
\begin{equation}  \label{eq:active_set_decomposition}
\x_{t} = \sum_{\vv \in \Coreset_{t} } \alpha^{(t)}_{\vv} \vv, \quad \text{and} \quad \y_{t} = \sum_{\vv \in \Vertices_{t}} \xi^{(t)}_{\vv} M\vv \;,
\end{equation}
When we choose to use the A\FW\ Alg.~\ref{alg:AFW} as inner loop algorithm, it can choose an \emph{away} direction to remove mass from ``bad'' atoms in the active set, i.e. to reduce~$\alpha^{(t)}_{\vv}$ for some~$\vv$ (see L\ref{algLine:drop_step} of Alg.~\ref{alg:AFW}), thereby avoiding the zig-zagging phenomenon that prevents \FW\ from achieving a faster convergence rate \citep{lacoste-julien_global_2015}.
On the other hand, the maximal step size for an \emph{away} step can be quite small ($\stepmax = \nicefrac{\alpha^{(t)}_{\vv}}{1-\alpha^{(t)}_{\vv}}$, where $\alpha^{(t)}_{\vv}$ is the weight of the away vertex in~\eqref{eq:active_set_decomposition}), yielding  to arbitrary small suboptimality progress when the line-search is truncated to such small step-sizes. A step removing an atom from the active set is called a \emph{drop step} (this is further discussed in Appendix~\ref{app:frank_wolfe_algorithms}), and Alg.~\ref{alg:AFW} loops until a non-drop step is obtained.
It is important to be able to upper bound the cumulative number of drop-steps in order to guarantee the termination of the inner loop  Alg.~\ref{alg:AFW} (Alg.~\ref{alg:AFW} ends only when it performs a non-drop step).
In App.~\ref{sub:upper_bound_on_the_number_of_drop_steps} we prove that the cumulative number of drop-steps after $t$ iterations cannot be larger than $t+1$.

\begin{figure}[t]
\centering
\fbox{\scalebox{.95}{
\begin{minipage}{.99\columnwidth}
\smallskip
\centerline{\bfseries FW Augmented Lagrangian method (\FWAL)}
\smallskip
At each iteration $t \ge 1$, we update the primal variable blocks $\x_t$ with a Frank-Wolfe step and then update the dual multiplier $\y_t$ using the updated primal variable:
\begin{equation}\label{eq:fwal}
\left\{\begin{array}{l}\displaystyle
\x_{t+1}=\mathcal{FW}(\x_{t};\LL(\cdot,\y_{t}))\;,\\[10pt]
\displaystyle \y_{t+1} = \y_{t} + \eta_{t}M \x_{t+1} \;,
\end{array}
\right.
\end{equation}
where $\eta_t>0$ is the step size for the dual update and $\mathcal{FW}$ is either Alg.~\ref{alg:AFW} or Alg.~\ref{alg:FW}  (see more in App.~\ref{app:frank_wolfe_algorithms}).
\end{minipage}
}}
\end{figure}
\section{Analysis of \FWAL{}} %
\label{sec:analysis_of_fw}

Solutions of \eqref{eq:opt_bis} are called saddle points, equivalently a vector $(\x^*,\y^*) \in \X \times \R^d$ is said to be a saddle point if the following is verified for all $(\x,\y) \in \X \times \R^d \;,$
\begin{equation}
  \LL(\x^*,\y) \leq \LL(\x^*,\y^*) \leq \LL(\x,\y^*)\;.
\end{equation}
Our assumptions (convexity of $f$ and $\X$, feasibility of $M\x = 0$, and crucially boundedness of $\X$) are sufficient for strong duality to hold~\citep[Exercise 5.25(e)]{boydBook}. Hence, the set of saddle points is not empty and is equal to $\X^* \times \Y^*$, where $\X^*$ is the set of minimizer of \eqref{eq:opt} and $\Y^*$ the set of maximizers of the augmented dual function $d$:
\begin{equation}\label{eq:def_dual_function}
d(\y) := \min_{\x \in \X} \LL(\x,\y)  \;.
\end{equation}
One of the issue of \AL{}\ is that it is a non-feasible method and thus  the function suboptimality is no longer a satisfying convergence criterion (since it can be negative). In the following section, we explore alternatives criteria to get a sufficient condition of convergence.
\subsection{Convergence Measures} %
\label{sub:convergence_measures}

Variants of \AL{} (also known as the methods of multipliers) update at each iteration both the primal variable and the dual variable. For the purpose of analyzing the popular \ADMM\ algorithm, \citet{hong_linear_2017} introduced two positive quantities which they called primal and dual gaps that we re-use in the analysis of our algorithm. Let $\x_t$ and $\y_t$ be the current primal and dual variables after $t$ iterations of the \FWAL\ algorithm~\eqref{eq:fwal}, the dual gap is defined as
\begin{equation}\label{eq:dd-gap}
\Delta_t^{(d)}:=d^*-d(\y_{t}) ~ \text{ where } \,d(\y_t) := \min_{\x \in \X} \LL(\x,\y_t)\;
\end{equation}
and $d^* := \max_{\y \in \R^d} d(\y)$. It represents the dual suboptimality at the $t$-th iteration. On the other hand, the ``primal gap'' at iteration $t$ is defined as
\begin{equation}\label{eq:pp-gap}
\Delta_t^{(p)} := \LL(\x_{t+1}, \y_{t})-d(\y_{t}),\quad  t\geq 0 \;.
\end{equation}
Notice that $\Delta_t^{(p)}$ is \emph{not} the suboptimality associated with the primal function $p(\cdot) := \max_{\y \in \R^d} \LL(\cdot,\y)$ (which is infinite for every non-feasible $\x$).
In this paper, we also define the shorthand,
\begin{equation}
\Delta_{t} := \Delta_t^{(p)}+\Delta_t^{(d)}\;.
\end{equation}
It is important to realize that since \AL\ is a non-feasible method, the standard convex minimization convergence certificates could become meaningless. In particular, the quantity $f(\x_{t}) - f^*$ might be negative since $\x_{t}$ does not necessarily belong to the constraint set of \eqref{eq:opt}. 
This is why it is important to consider the feasibility $\|M\x\|^2$.

In their work, \citet{yen_convex_2016,yen_dual_2016,huang_greedy_2017} only provide a rate on both gaps \eqref{eq:dd-gap} and \eqref{eq:pp-gap} which is not sufficient to derive guarantees on either how close an iterate is to the optimal set or how small is the suboptimality of the closest feasible point.
In this paper, we also prove the additional property that the feasibility $\|M\x\|^2$ converges to $0$ as fast as $\Delta_t$. 
But even with these quantities vanishing, the suboptimality of the closest feasible point can be significantly larger than the suboptimality of a point \mbox{$\epsilon$-close} to the optimum.
Concretely, let $\x \in \X$ and let $\tilde \x$ be its projection onto $\{\x \in \X \;|\; M\x = 0\}$, since $f$ is $L$-smooth we know that,
\begin{equation} \label{eq:conv_measure_smooth}
   |f(\tilde \x) - f(\x) - \innerProd{\nabla f(\x)}{\tilde \x- \x}| \leq \frac{L}{2} \|\x - \tilde \x\|^2 \,.
\end{equation}
On one hand, if the gradient is large and its angle with $\tilde \x-\x$ is not too small, $f(\tilde \x)$ may be significantly larger than $f(\x)$. On the other hand, if $\nabla f(\x)$ is not too large, we can upper bound the suboptimality at $\tilde \x$. Concretely, by \eqref{eq:conv_measure_smooth} we get,
\begin{equation}
  f(\tilde \x) \leq f(\x) + \|\nabla f(\x)\|\|\x-\tilde \x\| +  \frac{L}{2} \|\x - \tilde \x\|^2 \,.
\end{equation}
Moreover, since $\tilde \x$ is the projection of $\x$ onto the nullspace of $M$ we have that,
\begin{equation}
  \frac{\|M\x\|}{\sigma_{\max}(M)} \leq \|\x-\tilde \x\| \leq \frac{\|M\x\|}{\sigma_{\min}(M)} \,.
\end{equation}
Then, if $\|M\x\| \leq \epsilon$ and $f(\x) \leq \epsilon$ we have that 
\begin{equation}\label{eq:subopt_projection}
  f(\tilde \x) \leq (1+  \tfrac{\|\nabla f(\x)\|}{\sigma_{\min}(M)} + \tfrac{L \epsilon}{2\sigma_{\min}(M)^2}) \epsilon \,.
\end{equation}
The bound~\eqref{eq:subopt_projection} is not practical when the function appears to have gradients with large norms (which can be the case even close to the optimum for constrained optimization) or when $M$ appears to have small non-zero eigenvalues.
This is why we also consider the case where $f$ is strongly convex, allowing us to provide a bound on the distance to the optimum $\x^*$ (unique due to strong convexity).

\subsection{Properties of the augmented Lagrangian dual function} %
\label{sub:properties_of_the_augmented_dual_function}
The augmented dual function plays a key role in our convergence analysis.
One of our main technical contribution is the proof of a new property of this function which can be understood as a growth condition. This property is due to the smoothness of the objective function and the compactness of the constraint set. We will need an additional technical assumption called \emph{interior qualification} (a.k.a \emph{Slater's conditions}).
\begin{assumption}\label{assump:main}
$\exists \, \x^{(k)} \in \relint(\X_k), \; k \in [K]$ s.t. $\sum_{k=1}^K A_k \x^{(k)} =0$. 
\end{assumption}
Recall that $\x \in \relint(\X)$ if and only if $\x$ is an interior point relative to the affine hull of $\X$.
This assumption is pretty standard and weak in practice. It is a particular case of constraint qualifications~\citep{holmes1975geometric,gowda1990comparison}.
With this assumption, we can deduce a global property on the dual function that can be summarized as a quadratic growth condition on a ball of size $L_\lambda D^2$ and a linear growth condition outside of this ball. The optimization literature named such properties \emph{error bounds}~\citep{pang_error_1997}.

\begin{theorem}[Error bound]\label{thm:dual_function_PL}
 Let $d$ be the augmented dual function~\eqref{eq:def_dual_function}, if $f$ is a $L$-smooth convex function, 
  $\X$ a compact convex set and if Assump.~\ref{assump:main} holds, then there exist a constant $\alpha>0$ such that for all $\y \in \R^d$,
\begin{equation}\label{eq:main_lb_dual_directions}
  d^* - d(\y)
  \geq \frac{\alpha^2}{2}\min \left\{
        \frac{\dist(\y,\Y^*)^2}{L_\lambda D^2},\dist(\y,\Y^*)
  \right\} , \!\!
\end{equation}
where $D := \max_{(\x,\x') \in \X^2}\|\x - \x'\|$ is the diameter of $\X$ and $L_\lambda := L + \lambda \|M^\top \! M\|$.
 \end{theorem}
This theorem, proved in App.~\ref{sub:properties_of_the_dual_function_}, is crucial to our analysis. In our \emph{descent lemma}~\eqref{main_eq:fund_descent_lemma}, we want to relate the gap decrease to a quantity proportional to the gap. A consequence of~\eqref{eq:main_lb_dual_directions} is a lower bound 
 of interest:~\eqref{eq:main_d-gap}.

\paragraph{Issue in previous proofs.} %
\label{par:issue_in_previous_proofs_}

In previous work, \citet[Theorem 2]{yen_convex_2016} have a constant called $R_Y$ in the upper bound of $\Delta_t$ which may be infinite and lead to the trivial bound $\Delta_t \leq \infty$. Actually, $R_Y$ is an upper bound on the distance of the dual iterate$\;\y_t\;$to the optimal solution set $\Y^*$ of the augmented dual function. Since this quantity is not proven to be bounded, an element is missing in the convergence analysis.
In their convergence proof,
\citet{yen_dual_2016} and  \citet{huang_greedy_2017} use Lemma 3.1 from \citet{hong_linear_2012} (which also appears as Lemma 3.1 in the published version~\citep{hong_linear_2017}). This lemma states a result not holding for all $\y \in \R^d$ but instead for $(\y_t)_{t\in\N}$, which is the sequence of dual variables computed by the ADMM algorithm used in \citep{hong_linear_2017}. This sequence cannot be assimilated to the sequence of dual variables computed by the \GDMM\ algorithm since the update rule for the primal variables in each algorithm is different: the primal variable are updated with \FW\ steps in one algorithm and with a proximal step in the other.
The properties of this proximal step are intrinsically different from the \FW\ steps computing the updates on the primal variables of \FWAL. To our knowledge, there is no easy fix (details in App.~\ref{sec:discussion_on_previous_proofs}) to get a similar result as the one claimed in \citep[Lem.~4]{yen_dual_2016} and \citep[Lem.~4]{huang_greedy_2017}.
\subsection{Specific analysis for \FWAL} %
\label{sub:analysis_for_fwal}
\paragraph{Convergence over general convex sets.} %
\label{par:general_convex_sets_}
When $\X$ is a general convex compact set and $f$ is $L$-smooth,
 Algorithms~\ref{alg:AFW} and \ref{alg:FW} are able to perform a decrease on the objective value proportional to the square of the suboptimality~\citep[Lemma 5]{jaggi_revisiting_2013}, \citep[(31)]{lacoste-julien_global_2015}, we will call this a \emph{sublinear decrease} since it leads to a sublinear rate for the suboptimality: for any $\x \in \mathcal \X, \, \y \in \R^d$  they compute $\x^+ := \mathcal{FW}(\x;\LL(\cdot,\y))$, such that for all $\gamma\in[0,1],$
\begin{equation}\label{eq:def_subl_decrease}
    \LL(\x^+ \! , \y) - \LL(\x,\y) \leq \gamma \!\left( d(\y) - \LL(\x,\y)  \right)  + \mfrac{\gamma^2L_\lambda D^2}{2},\!\!\!
\end{equation}
where $L_\lambda$ is the Lipschitz constant of $\nabla \LL$ and $D$ the diameter of $\X$.
Recall that $d(\y) := \min_{\x' \in \X}\LL(\x',\y)$.
Note that setting $\gamma = 0$ gives $\LL(\x^+,\y) \leq \LL(\x,\y)$ and optimizing the RHS respect to $\gamma$ yields a decrease proportional to $(d(\y) - \LL(\x,\y))^2$.
The \GDMM\ algorithm of \citet{yen_convex_2016,yen_dual_2016,huang_greedy_2017} relies on the assumption of $\X$ being polytope, hence we obtain under this general assumption of sublinear decrease a new result on \AL\ with \FW.
This result covers the case of the simultaneously sparse and low rank matrices~\eqref{eq:constrained_sparse_low_rank} where the trace norm ball is not a polytope.

\begin{theorem}[Rate of \FWAL{} with Alg.~\ref{alg:FW}]
\label{thm:main_general}
 Under Assumption~\ref{assump:main}, if $\X$ is a convex compact set and $f$ is a $L$-smooth convex function and $M$ has the form described in~\eqref{eq:form_M}, 
then using any algorithm with \emph{sublinear decrease}~\eqref{eq:def_subl_decrease} as inner loop in \FWAL\ \eqref{eq:fwal} and $\eta_t :=  \min\Big\{\frac{2}{\lambda}, \frac{\alpha^2}{2 \delta}\Big\}\frac{2}{t+2}$, we have that there exists a bounded $t_0 \geq 0$ such that $\forall t \geq t_1\geq t_0$,
  \begin{equation}
    \Delta_{t} \leq \frac{4\delta(t_0+2)}{t+2}, \,
    \min_{t_1 \leq s-1 \leq t} \! \|M\x_{s}\|^2
   \leq \frac{\bigO{}(1)}{t-t_1+1}
\end{equation}
where $D := \max_{\x, \x' \in \X}\|\x-\x'\|$ is the diameter of $\X$, $L_\lambda := L + \lambda \|M^\top \! M\|$ the Lipschitz constant of $\nabla \LL$, $\delta := L_\lambda D^2$ and $\alpha$ is defined in Thm.~\ref{thm:dual_function_PL}.
\end{theorem}
In App.~\ref{sub:proof_of_1_in_theorem_thm:main}, we provide an analysis for different step size schemes and explicit bounds on $t_0$.

\paragraph{Convergence over Polytopes.} %
\label{par:polytopes_}
On the other hand, if $\X$ is a polytope and $f$ a generalized strongly convex function, recent advances on \FW\ proposed global linear convergence rates using \FW\ with away steps~\citep{lacoste-julien_global_2015,garber_linear-memory_2016-1}.
Note that in the augmented formulation, $\lambda>0$ and thus $\frac{1}{2}\|M\cdot\|^2$ is a generalized strongly convex function, making $\LL(\cdot,\y)$ a generalized strongly convex function for any $\y \in \R^d$ (see App.~\ref{sub:constants_for_the_sublinear_and_geometric_decrease} for details).
We can then use such linearly convergent algorithms to improve the rate of \FWAL.
More precisely, we will use the fact that Algorithm~\ref{alg:AFW} performs a \emph{geometric decrease}~\citep[Theorem 1]{lacoste-julien_global_2015}: for $\x^+ := \mathcal{FW}(\x;\LL(\cdot,\y))$, there exists $ \rho_A< 1$ such that for all $\x \in \mathcal \X$ and $\y \in \R^d$,
  \begin{equation}\label{eq:def_geom_decrease}
   \LL(\x^+\!, \y) -  \LL(\x,\y) \leq  \rho_A  \big[\min_{\x' \in \X}\LL(\x',\y)  - \LL(\x,\y) \big]. \!\!\!
  \end{equation}
The constant $\rho_A$~\citep{lacoste-julien_global_2015} depends on the smoothness, the generalized strong convexity of $\LL(\cdot,\y)$ (does not depend on~$\y$, but depends on~$M$) and the condition number of the set $\X$ depending on its geometry (more details in App.~\ref{sub:constants_for_the_sublinear_and_geometric_decrease}).
\begin{theorem}[Rate of \FWAL{} with inner loop Alg.~\ref{alg:AFW}]
\label{thm:main_polytope}
 Under the same assumptions as in Thm.~\ref{thm:main_general} and if moreover $\X$ is a polytope and $f$ a generalized strongly convex function,
 then using Alg~\ref{alg:AFW} as inner loop and a constant step size $\eta_t = \frac{\lambda \rho_A}{4}$, the quantity $\Delta_t$ decreases by a uniform amount for finite number of steps $t_0$ as,
\begin{equation}
    \Delta_{t+1}-\Delta_{t}
  \leq  -\frac{\lambda \alpha^2 \rho_A}{8}\;,
\end{equation} until $\Delta_{t_0} \leq L_\lambda D^2$. Then for all $t\geq t_0$ we have that the gap and the feasibility violation decrease linearly as,
  \begin{equation*}
    \Delta_{t} \leq \frac{\Delta_{t_0}}{(1+\rho)^{t-t_0}}\,,
     \; \; \|M\x_{t+1}\|^2 \leq  \frac{16}{\lambda \cdot \rho_A} \frac{\Delta_{t_0}}{(1+\rho)^{t-t_0}}\;,
  \end{equation*}
  where $\rho : =  \min \big\{ \frac{\rho_A}{2}, \frac{\rho_A\lambda \alpha^2}{8 L_\lambda D^2}\big\}$ and $L_\lambda := L + \lambda\|M^\top \!M\|$.

\end{theorem}

\paragraph{Strongly convex functions.} %
\label{par:strongly_convex_functions_}
When the objective function $f$ is strongly convex, we are able to give a convergence rate for the distance of the primal iterate to the optimum. As argued in Sec.~\ref{sub:convergence_measures}, an iterate close to the optimal point lead to a ``better'' approximate solution than an iterate achieving a small gap value. 
\begin{theorem}\label{thm:strong}
Under the same assumptions as in Thm.~\ref{thm:main_general}, if $f$ is a $\mu$-strongly convex function, then the set of optimal solutions $\X^*$ is reduced to $\{\x^*\}$ and for any $t \geq t_1 \geq 8t_0+14$,
   \begin{equation}
     \min_{t_1 +1 \leq s \leq t +1} \|\x_t - \x^*\|^2 \leq \frac{4}{\mu} \frac{\bigO{}(1)}{t-t_1 +1}\;.
   \end{equation}
Moreover, if $\X$ is a compact polytope, and if we use Alg.~\ref{alg:AFW}, then the distance of the current point to the optimal set vanishes as (with $\rho$ as defined in Thm.~\ref{thm:main_polytope}):
  \begin{equation}
    \|\x_{t+1} -\x^*\|^2  \leq \frac{2\Delta_{t_0}(\sqrt{2}+1)}{\mu(\sqrt{1+\rho})^{t-t_0}} +\frac{O(1)}{(1+\rho)^{t-t_0}} \;.
  \end{equation}
\end{theorem}
This theorem is proved in App.~\ref{app:proof_of_thm} (Cor.~\ref{cor:general_polytope_feasibility} and Cor.~\ref{cor:main_polytope}).
For an intersection of sets, the three theorems above give stronger results than~\citep{yen_dual_2016,huang_greedy_2017}
since we prove that the distance to the optimal point as well as the feasibility vanish linearly.

\proof[\textbf{\emph{Proof sketch of Thm~\ref{thm:main_general} and~\ref{thm:main_polytope}}}]

Our goal is to obtain a convergence rate on the sum gaps \eqref{eq:dd-gap} and \eqref{eq:pp-gap}.
First we show that the dual gap verifies
\begin{equation}
\Delta_{t+1}^{(d)}-\Delta_{t}^{(d)}\leq -\eta_{t} \innerProd{M\x_{t+1}}{M\hat \x_{t+1}} \; \label{main_eq:dual-gap}
\end{equation}
where $\hat \x_{t+1} := \argmin_{\x\in \X} \LL(\x,\y_{t+1})$. Similarly, we prove the following inequality for the primal gap
\begin{align}
\label{main_eq:primal-gap}
\Delta_{t+1}^{(p)}-\Delta_{t}^{(p)}
&\leq \eta_{t}\|M\x_{t+1}\|^2 \nonumber\\
& \quad\,+ \LL(\x_{t+2},\y_{t+1})-\LL(\x_{t+1},\y_{t+1}) \nonumber\\
& \quad\,-\eta_{t} \innerProd{M\x_{t+1}}{M\hat \x_{t+1}}\;.
\end{align}

Summing \eqref{main_eq:dual-gap} and \eqref{main_eq:primal-gap} and using that $\|M\x_{t+1}-M \hat\x_{t+1}\|^2 \leq  \frac{2}{\lambda} \left( \LL(\x_{t+1},\y_{t+1}) - \LL(\hat \x_{t+1},\y_{t+1}) \right)$, we get the following \emph{fundamental descent lemma},
  \begin{align}
   \Delta_{t+1}-\Delta_{t} \notag
   &\leq \LL(\x_{t+2},\y_{t+1})-\LL(\x_{t+1},\y_{t+1}) \\
   & \quad \, + \frac{2 \eta_{t}}{\lambda} \left( \LL(\x_{t+1},\y_{t+1})-\LL(\hat\x_{t+1},\y_{t+1}) \right) \notag \\
   & \quad \, - \eta_{t}\|M\hat \x_{t+1}\|^2\;. \label{main_eq:fund_descent_lemma}
  \end{align}
 We now crucially combine \eqref{eq:main_lb_dual_directions} in Thm.~\ref{thm:dual_function_PL} and the fact that $\Delta^{(d)}_t \leq \dist(\y^t,\Y^*) \|M\hat \x_{t+1}\|$ to obtain,
\begin{equation}\label{eq:main_d-gap}
  \tfrac{\alpha^2}{2L_\lambda D^2}\min\{\Delta_{t+1}^{(d)},L_\lambda D^2 \} \leq  \|M\hat \x_{t+1}\|^2 \;,
\end{equation}
and then,
\begin{align}\label{eq:main_fun_descent}
   \Delta_{t+1}-\Delta_{t}
   &\leq \LL(\x_{t+2},\y_{t+1})-\LL(\x_{t+1},\y_{t+1}) \nonumber\\
   & \quad \, + \frac{2 \eta_{t}}{\lambda} \left( \LL(\x_{t+1},\y_{t+1})-\LL(\hat\x_{t+1},\y_{t+1}) \right) \nonumber\\
   & \quad \, - \eta_{t} \tfrac{\alpha^2}{2L_\lambda D^2}\min\{\Delta_{t+1}^{(d)},L_\lambda D^2 \} \;.
\end{align}
Now the choice of the algorithm to get $\x_{t+2}$ from $\x_{t+1}$ and $\y_{t+1}$ is decisive:

If $\X$ is a polytope and if an algorithm with a \emph{geometric decrease}~\eqref{eq:def_geom_decrease} is used, setting $\eta_{t} = \frac{\lambda \cdot \rho_A}{4}$ we obtain
\begin{equation}\label{main_eq:fund descent} \notag
\begin{aligned}
  \Delta_{t+1}-\Delta_{t}
  & \leq -\frac{\rho_A}{2 } \left( \LL(\x_{t+1},\y_{t+1}) - \LL(\hat \x_{t+1},\y_{t+1}) \right) \\
  & \quad\, -\frac{\lambda \cdot \rho_A}{4}\|M \x_{t+1}\|^2 \;.
\end{aligned}
\end{equation}
Since $\LL(\x_{t+2},\y_{t+1}) \leq  \LL(\x_{t+1},\y_{t+1})$ (L\ref{line:line_search}), we have
\begin{equation}
\Delta_{t+1}^{(p)} \leq \LL(\x_{t+1},\y_{t+1}) - \LL(\hat \x_{t+1},\y_{t+1}) \;,
\end{equation}
leading us to a geometric decrease for all $t \geq t_0$,
\begin{equation}
  \Delta_{t+1} \leq \frac{\Delta_{t}}{1+\rho} \;\; \text{where} \;\; \rho :=  \tfrac{\rho_A}{2} \min \left\{ 1, \tfrac{\lambda \alpha^2}{8 L_\lambda D^2}\right\}.
\end{equation}
Additionally we can deduce from \eqref{main_eq:fund_descent_lemma} that,
  \begin{equation}
    \eta_{t}\|M\hat\x_{t+1}\|^2 \leq \Delta_{t}
    \quad \text{and} \quad
    \eta_{t}\|M\x_{t+1}\|^2  \leq 4\Delta_{t} \;.
  \end{equation}
\begin{figure*}[h]
    \centering
    \begin{subfigure}[b]{0.33\linewidth}
    \hspace*{-3mm}
        \centering
        \includegraphics[width = 1.1 \linewidth]{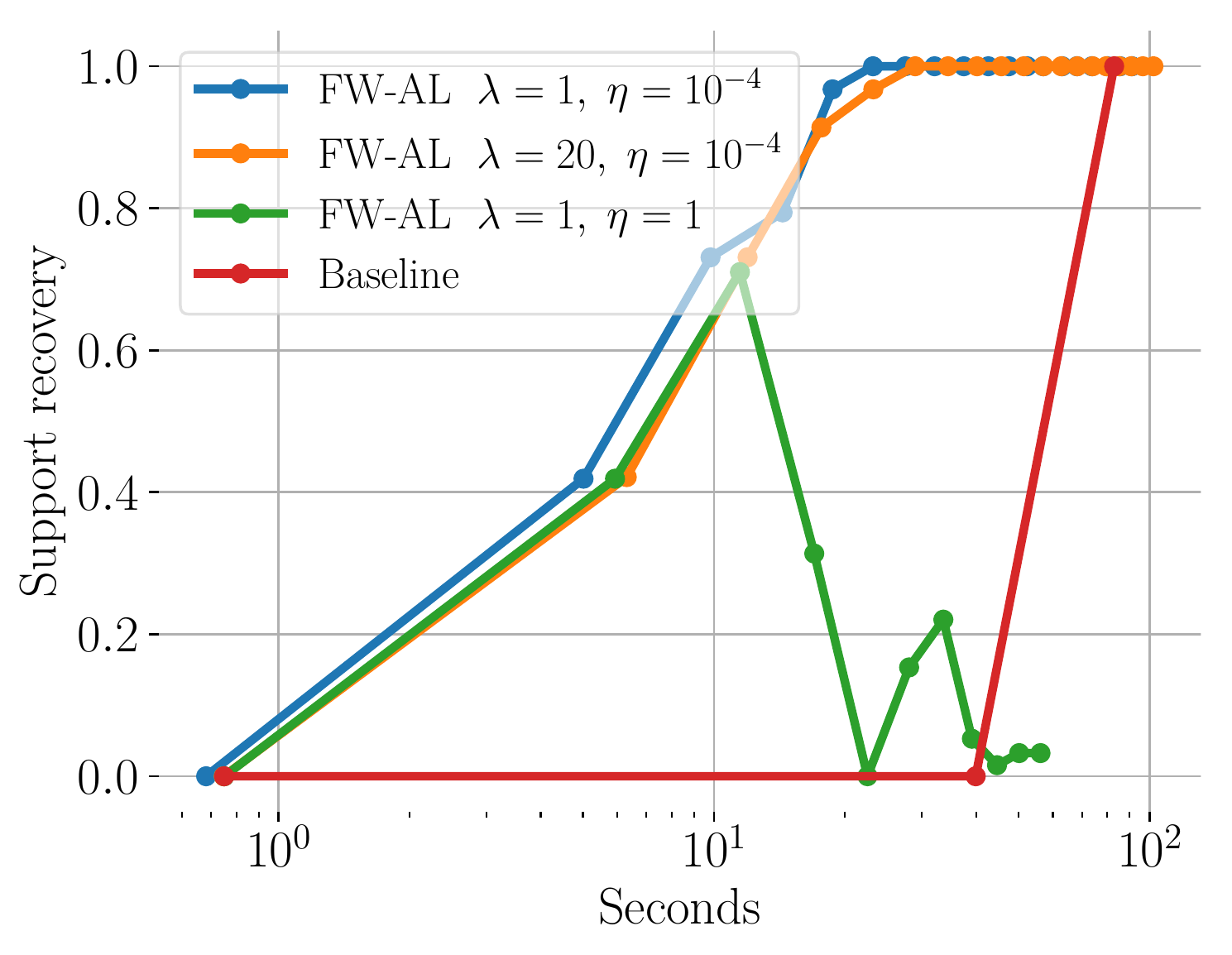}
        \caption{\small Fraction of the support recovered.}
        \label{fig:support_recovery}
    \end{subfigure}
    \;\; \;
    \begin{subfigure}[b]{0.30\linewidth}
        \centering
        \includegraphics[width = 1 \linewidth]{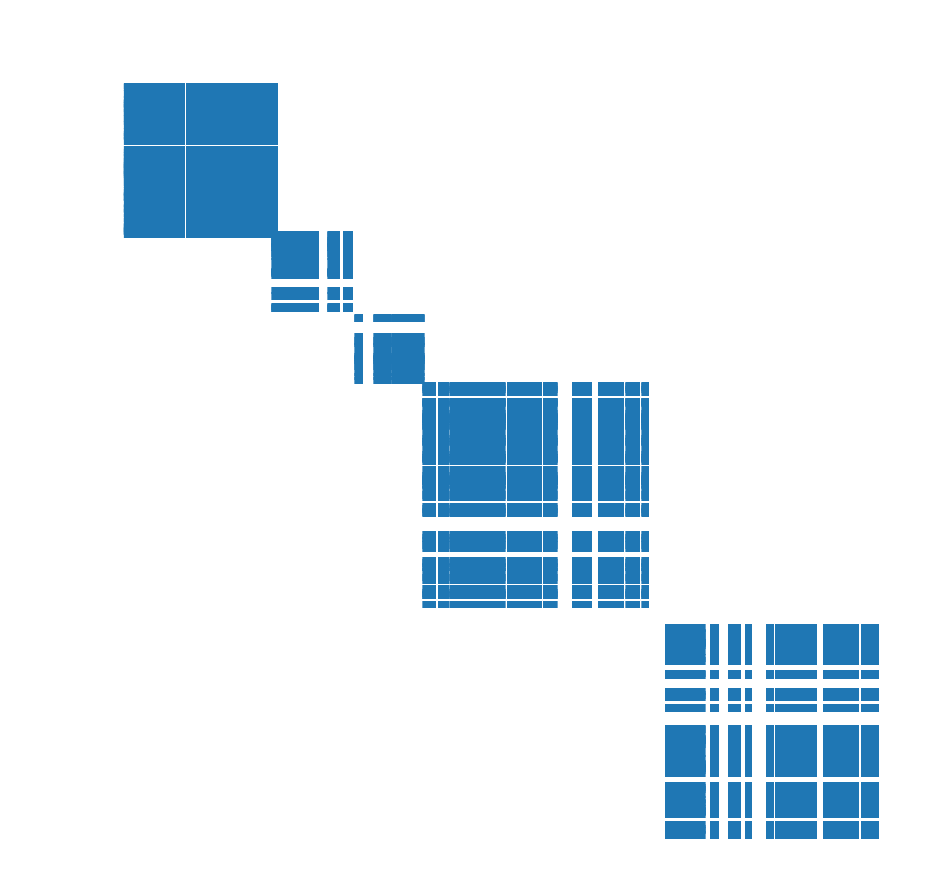}
        \caption{ \small Matrix recovered with \FWAL.}
        \label{fig:fw_recover}
    \end{subfigure}
    \; \;
    \begin{subfigure}[b]{0.30\linewidth}
        \centering
        \includegraphics[width = 1 \linewidth]{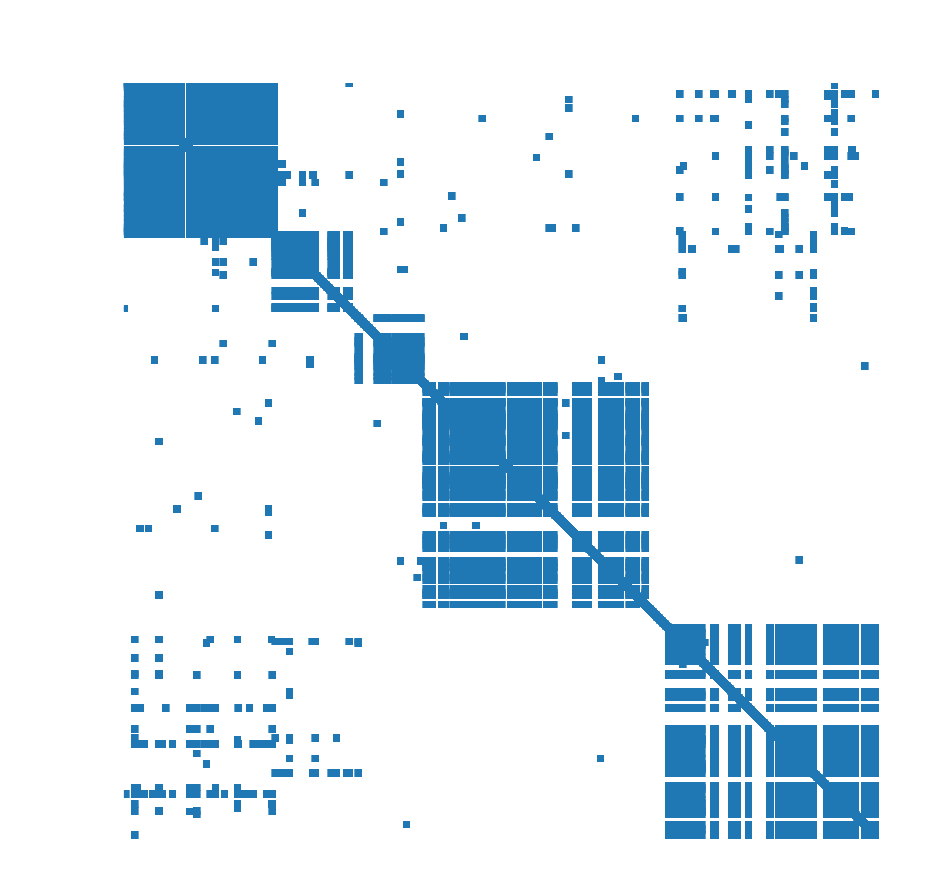}
        \caption{ \small Matrix recov. with the baseline.}
        \label{fig:split_recover}
    \end{subfigure}
    \caption{\small Fig.~\ref{fig:support_recovery} represent the fraction of the support of $\Sigma$ recovered as a function of time ($d^2 = 1.6 \cdot 10^{7}$ and the matrix computed is thresholded at $10^{-2}$). The baseline is the generalized forward backward algorithm. \FWAL{} requires a small enough step size $\eta$ to recover the support otherwise it diverges (green curve) and does not require a lot of tuning for $\lambda$ (blue and orange curve).
    Fig~\ref{fig:fw_recover} and \ref{fig:split_recover} compare the matrices recovered for $d^2 = 10^6$ after one minute of computation.
    }\label{fig:experiments_support_recovery}
\end{figure*}
  If $\X$ is not a polytope, we can use an algorithm with a \emph{sublinear decrease}~\eqref{eq:def_subl_decrease} to get from \eqref{eq:main_fun_descent} that $\forall t\geq 0 \;,$
  \begin{equation}
      \Delta_{t+1} - \Delta_{t} \leq - a \eta_{t}\min\{\Delta_{t+1},\delta\} + (a\eta_{t})^2\tfrac{C}{2}\;,
  \end{equation}
  where $a,\delta$ and $C$ are three positive constants.
  Setting $\eta_{t} = \tfrac{2}{a(t+2)}$ we can prove that there exists $t_0 \geq \frac{C}{\delta}$ s.t.,

  \begin{equation}
    \Delta_{t+1} \leq \frac{4 \delta (2 + t_0) }{(t+2)}, \quad \forall t\geq t_0 \;.
  \end{equation}
  Providing the claimed convergence results.
\endproof
\section{Illustrative Experiments} %
\label{sec:experiments}
\begin{figure}
\centering
       \includegraphics[width = .87 \linewidth]{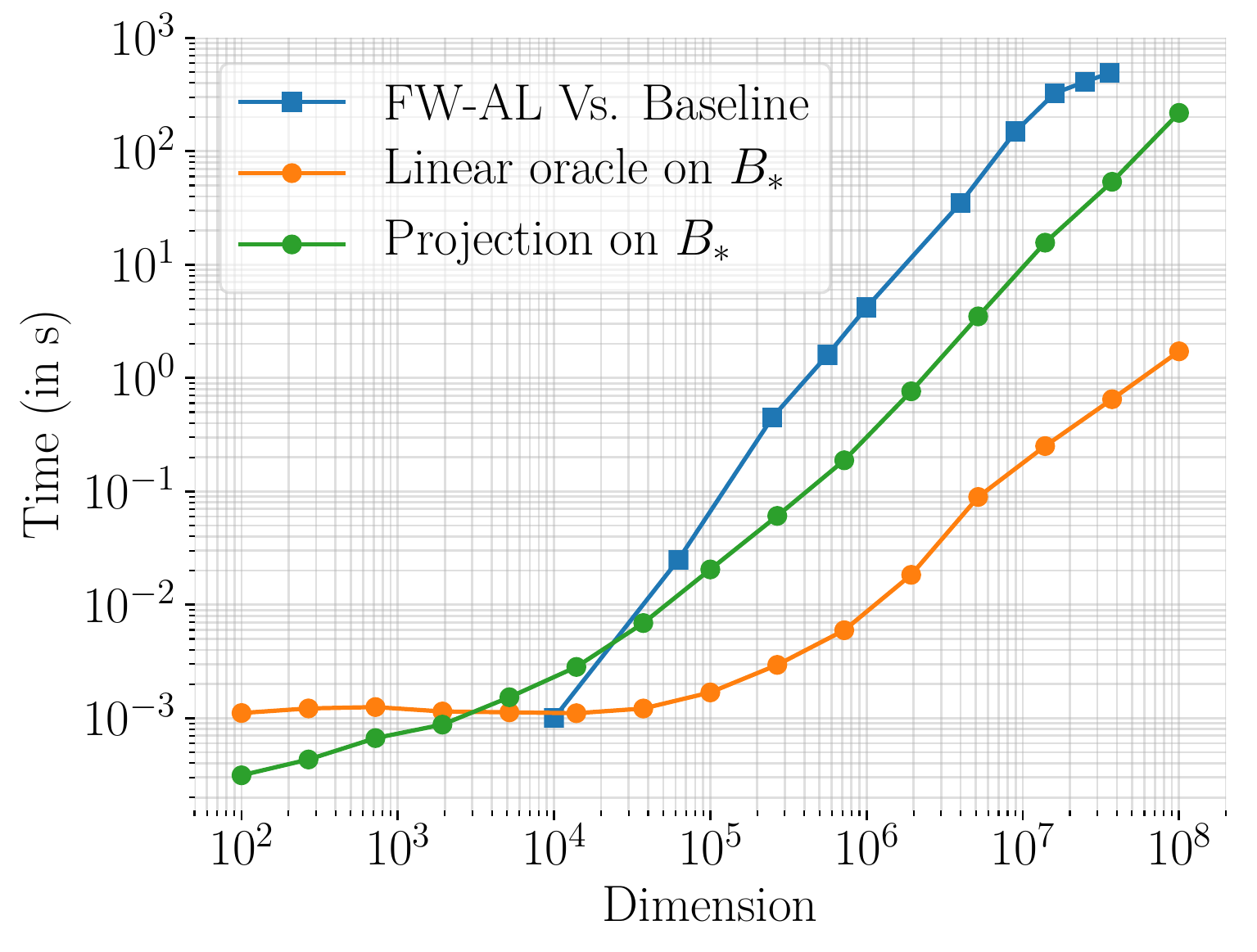}
        \caption{ \small Time complexity of the LMO vs. the projection on the trace norm ball. The blue curve represents the time spent by the generalized forward backward algorithm to reach a better point than the one computed by \FWAL.}
        \label{fig:experiments_time}
\end{figure}
Recovering a matrix that is simultaneously low rank and sparse has applications in problems such as covariance matrix estimation, graph denoising and link prediction~\citep{richard_estimation_2012}.
We compared \FWAL\ with the proximal splitting method on a covariance matrix estimation problem.
We define the $\|\cdot \|_1$ norm of a matrix $S$ as $\|S\|_1 := \sum_{i,j} |S_{i,j}|$ and its trace norm as $\|S\|_* := \sum_{i=1}^{\text{rank}(S)} \sigma_i$, where $\sigma_i$ are the singular values of $S$ in decreasing order. Given a symmetric positive definite matrix $\hat\Sigma \in \R^{d\times d}$ we use the square loss as strongly convex objective for our optimization problem,
  \begin{equation}\label{eq:constrained_sparse_low_rank}
  \min_{S \succeq 0, \|S\|_1 \leq \beta_1 , \|S\|_{*} \leq \beta_2 }  \;\|S-\hat\Sigma\|_2^2\;.
\end{equation}
The linear oracle for $\X_1 := \{S \succeq 0,\, \|S\|_1 \leq \beta_1\}$ is
\begin{equation} \notag
\LMO_{\X_1}(D) := \beta_1\tfrac{E_{ij} + E_{ji}}{2}, \; (i,j) \in \argmin_{(i,j) \in [d]\times[d]} D_{i,j}+D_{j,i}
\end{equation}
 where $(E_{ij})$ is the standard basis of $\R^{d \times d}$. The linear oracle for $\X_2: = \{S\succeq 0,\,\|S\|_* \leq \beta_2\}$ is
\begin{equation}\label{eq:LMO_trance_norm}
\LMO_{\X_2}(D) := \beta_2 \cdot  U_1^\top U_1\; ,
\end{equation}
where $D = \left[ U_1,\ldots,U_d \right]  \diag(\sigma_1,\ldots,\sigma_d) \left[ U_1,\ldots,U_d \right] ^\top\!.$ For this problem, the matrix $D$ is always symmetric because the primal and dual iterates are symmetric as well as the gradients of the objective function.
Eq.~\eqref{eq:LMO_trance_norm} can be computed efficiently by the Lanczos algorithm~\citep{paige1971computation, kuczynski_estimating_1992-1} whereas the forward backward splitting which is the standard splitting method to solve~\eqref{eq:constrained_sparse_low_rank} needs to compute projections over the trace norm ball via a complete diagonalization which is $\bigO(d^3)$. For large $d$, the full diagonalization becomes untractable, while the Lanczos algorithm is more scalable and requires less storage (see Fig.~\ref{fig:experiments_time}).

The experimental setting is done following~\citet{richard_estimation_2012}: 
we generated a block diagonal covariance matrix~$\Sigma$ to draw $n$ vectors $\x_i \sim \mathcal N(0,\Sigma)$. We use $5$ blocks of the
form $\vv \vv^\top$ where $\vv \sim \mathcal U ([-1,1])$.
In order to enforce sparsity, we only kept the entries $(i,j)$ such that $|\Sigma_{i,j}| > .9$.
Finally, we add a gaussian noise $\mathcal N(0,\sigma)$ on each entry $\x_i$ and observe $\hat{\Sigma}=\sum^n_{i=1}\x_i \x_i^T$. In our experiment $n = d, \sigma = 0.6$.
We apply our method, as well as the the generalized forward backward splitting used by \citet{richard_estimation_2012}. This algorithm is the baseline in our experiments. It has been originally introduced by~\citet{raguet2013generalized}, to optimize \eqref{eq:constrained_sparse_low_rank}
performing projections over the constraint sets.
The results are presented in Fig.~\ref{fig:experiments_support_recovery} and~\ref{fig:experiments_time}.
We can say that our algorithm performs better than the baseline for high dimensional problems for two reasons: in high dimensions, only one projection on the trace norm ball $B_{*}$ can take hours (green curve in Fig.~\ref{fig:experiments_time}) whereas solving a \LMO\ over $B_{*}$ takes few seconds. Moreover, the iterates computed by \FWAL\ are naturally sparse and low rank, so we then get a better estimation of the covariance matrix at the beginning of the optimization as illustrated in Fig.~\ref{fig:fw_recover} and \ref{fig:split_recover}.

\subsection*{Acknowledgments}
We thank an anonymous reviewer for valuable comments which enabled us to improve the proofs.
This research was partially supported by the Canada Excellence Research Chair in “Data Science for Realtime Decision-making”, by the NSERC Discovery Grant RGPIN-2017-06936 and by the European Union's Horizon 2020 research and innovation program under the Marie Sk{\l}odorowska-Curie grant agreement 748900.

\bibliography{Splitting}
\bibliographystyle{abbrvnat}
\newpage
\appendix
\onecolumn

\section{Frank-Wolfe inner Algorithms} %
\label{app:frank_wolfe_algorithms}

\subsection{Upper bound on the number of drop-steps} %
\label{sub:upper_bound_on_the_number_of_drop_steps}

\begin{proposition}[Sparsity of the iterates and upper bound on the number of drop-steps]\label{prop:sparse_drop_step} The iterates computed by \FWAL\ have the following properties,
\begin{enumerate}
  \item After $t$ iterations, the iterates $\x_{t}$ (resp. $\y_{t}$) are a convex (resp. conic) combination of their initialization and the oracle's outputs (resp. times $M$) for the first $t$ iterations.
  \item
  If the algorithm $\mathcal{FW}$ is A\FW\ (Alg.~\ref{alg:AFW}), and if we initialize our algorithm at a vertex, after $t$ iterations of the main loop the cumulative number of drop-steps performed in the inner algorithm~\ref{alg:AFW} is upper bounded by $t+1$.
\end{enumerate} \label{thm:sparsity}
\end{proposition}
\proof
The first point comes from~\eqref{eq:active_set_decomposition}.

A \emph{drop step} happens when $\gamma_t=\stepmax$ in the away-step update L.~\ref{line:line_search} of Alg.~\ref{alg:AFW}. In that case, at least one vertex is removed from the active set.
The upper bound on the number of drop step can be proven with the same technique as in \citep[Proof of Thm.~8]{lacoste-julien_global_2015}.
Let us call $A_t$ the number of \FW\ steps (which potentially adds an atom in $\Coreset_t$) and $D_t$ the number of \emph{drop-steps}, i.e., the number of \emph{away steps} where at least one atom from $\Coreset_t$ have been removed (and thus $\stepsize_t = \stepmax$ for these). Considering \FWAL\ with A\FW\, after $t$ iterations we have performed $t$ non drop-steps in the inner loop, since it is the condition to end the inner loop, then
\begin{equation} \label{eq:NumberSteps}
A_t \leq t, \quad \text{and}\quad A_t -  D_t + |\Coreset_0| \geq |\Coreset_t| \geq 0 \;.
\end{equation}
Since by assumption $|\Coreset_0|=1$, this leads directly to $D_t \leq A_t +1 \leq t+1.$
\endproof

\subsection{Other \FW\ Algorithms Available} %
\label{sub:other_fw_algorithm_available}
Any Frank-Wolfe algorithm performing \emph{geometric decrease}~\eqref{eq:def_geom_decrease} or \emph{sublinear decrease}~\eqref{eq:def_subl_decrease} can be used an inner loop algorithm. For instance, the block-coordinate Frank-Wolfe method~\citep{lacostejulien:hal-00720158} performs a sublinear decrease in expectation and the fully-corrective
Frank-Wolfe method~\citep{lacoste-julien_global_2015} or \citet{garber_linear-memory_2016-1}'s algorithm perform a geometric decrease.

\subsection{Constants for the sublinear and geometric decrease} %
\label{sub:constants_for_the_sublinear_and_geometric_decrease}
In order to be self-contained, we will introduce the definitions of the constants introduced in the definition of \emph{sublinear decrease}~\eqref{eq:def_subl_decrease} and \emph{geometric decrease}~\eqref{eq:def_geom_decrease}.
\paragraph{Sublinear Decrease.} %
\label{par:sublinear_decrease_}
Let us first recall Equation~\eqref{eq:def_subl_decrease} describing the sublinear decrease:
\begin{equation}\notag
    \LL(\x^+ , \y) - \LL(\x,\y) \leq -\gamma \left( \LL(\x,\y) -  \min_{\x' \in \X}\LL(\x',\y)  \right)  +\gamma^2 \frac{L_\lambda D^2}{2} \;.
\end{equation}
The sublinear decrease is a consequence of the standard \emph{descent lemma}~\citep[(1.2.5)]{nesterov_introductory_2004}.
The constant $L_\lambda$ is the smoothness of $\LL$ and $D$ the diameter of $\X$.
This property has been proved for the block-coordinate Frank-Wolfe method\footnote{For BCFW, the sublinear decrease is valid on the expectation of the suboptimality, then the proofs with this algorithm as an inner-loop require a bit of extra work.}~\citep{lacostejulien:hal-00720158}, usual Frank-Wolfe~\citep{jaggi_revisiting_2013} and Frank-Wolfe with away-step~\citep{lacoste-julien_global_2015}.

If $f$ is $L$-smooth we have that the function $\LL(\cdot,\y)$ is $L_\lambda :=L + \lambda \|M^\top \! M\|$-smooth for any $\y \in \R^d$, and then,
\begin{equation}
  L_\lambda D^2 \leq \big(L + \lambda\|M^\top \!M \|\big) D_\X^2 \;.
\end{equation}
Recall that, for matrices $\|\cdot\|$ is the spectral norm.

\paragraph{Geometric Decrease.} %
\label{par:geometric_decrease_}
If the function $f$ is a generalized strongly convex function, then $\LL(\cdot,\y)$ is also a generalized strongly convex function. More generally, let $h_1$ and $h_2$ be two generalized strongly convex functions. Then according to the definition~\eqref{eq:generalized_strong_conv}, there exist $A_1,A_2,b_1,b_2$ and two strongly convex functions $g_1,g_2$ such that, $h_1(\x) =  g_1(A_1\x) + \innerProd{\bm{b}_1}{\x}$ and $h_2(\x) = g_2(A_2\x) + \innerProd{\bm{b}_2}{\x}$. Thus,
\begin{equation}
  h_1(\x) +h_2(\x) 
  = g_1(A_1\x) + \innerProd{\bm{b}_1}{\x} + g_2(A_2\x) + \innerProd{\bm{b}_2}{\x}
  = g(A\x) + \innerProd{\bm{b}}{\x} \,,
\end{equation}
where $A\x = [A_1\x;A_2\x],\; \bm{b} = [\bm{b}_1;\bm{b}_2]$ and $g([\uu;\vv]) = g_1(\uu) + g_2(\vv)$. The function $g$ is strongly convex by strong convexity of $g_1$ and $g_2$.

We can say that since $\LL(\cdot,\y)$ is a generalized strongly convex function (with a constant uniform on $\y$) and $\X$ a polytope, we have the geometric descent lemma from \citet[Theorem 1]{lacoste-julien_global_2015}. The constant $\rho_A$ is the following
\begin{equation}
  \rho_A := \frac{\mu_\lambda}{4L_\lambda} \left(\frac{\delta_\X}{D_\X}\right)^2 \;,
\end{equation}
where $\mu_\lambda$ and $L_\lambda$ are respectively the generalized strong convexity constant~\citep[Lemma 9]{lacoste-julien_global_2015} and the smoothness constant of $\x \mapsto  f(\x) + \tfrac{\lambda}{2}\|M\x\|^2$, and $D_\X$ and $\delta_\X$ are respectively the diameter and the pyramidal width of $\X$ as defined by \citet{lacoste-julien_global_2015}. Note that if $M$ is full rank, the strong convexity constant $\mu$ is lower bounded by $\lambda \sigma_{min}^2$ where $\sigma_{min}^2$ is the smallest singular value of $M$. Otherwise, if $M$ is not full rank, one can still use the lower bound on the generalized strong convexity constant given by \citet[Lemma 9]{lacoste-julien_global_2015}.

\section{Previous work} %
\label{sec:previous_work}

\subsection{Discussion on previous proofs} %
\label{sec:discussion_on_previous_proofs}

The convergence result stated by \citet[Theorem 2]{yen_convex_2016} is the following (with our notation)
\begin{equation}
  \Delta^{(p)}_t + \Delta^{(d)}_t \leq \frac{\omega}{t}
  \quad \text{where} \quad
  \omega := \frac{4}{1-\rho_A}\max \left(  \Delta^{(p)}_0 + \Delta^{(d)}_0, 2 R_Y^2/\lambda \right) \;,
\end{equation}
and $R_Y := \sup_{t\geq 0} \dist(\y_t,\Y^*)$. This quantity was introduced in the last lines of the appendix without any mention to its boundedness. In our opinion, it is as challenging to prove that this quantity is bounded as to prove that $\Delta_t$ converges.

In more recent work, \citet{yen_dual_2016} and  \citet{huang_greedy_2017} use a different proof technique in order to prove a linear convergence rate for their algorithm.
In order to avoid getting the same problematic quantity $R_Y$, they use Lemma 3.1 from~\citep{hong_linear_2012} (which also appears as Lemma 3.1 in the published version~\citep{hong_linear_2017}). This lemma states a result not holding for all $\y \in \R^d$ but instead for $(\y_t)_{t\in\N}$, which is the sequence of dual variables computed by the algorithm introduced in \citep{hong_linear_2017}. This sequence cannot be assimilated to the sequence of dual variables computed by the \GDMM\ algorithm since the update rule for the primal variables in each algorithm is different: the primal variable are updated with \FW\ steps in one algorithm and with a proximal step in the other.
The properties of this proximal step are intrinsically different from the \FW\ steps computing the updates on the primal variables of \FWAL.
One way to adapt this Lemma for \FWAL{} (or \GDMM{}) would be to use \citep[Lemma 2.3 c]{hong_linear_2017}. Unfortunately, this result is local (only true for all $\y \in \Y$ such that $\|\nabla d(\y)\|\leq \delta$ with $\delta$ fixed), whereas a global result (true for all $\delta$) seems to be required with the proof technique used in \citep{yen_dual_2016,huang_greedy_2017}. It is also mentioned in  \citep[proof of Lemma 2.3 c]{hong_linear_2017} that ``if in addition $\y$ also lies in some compact set $\Y$, then the dual error bound hold true for all $\y \in \Y$'' then showing that $R_Y$ is bounded would fix the issue, but as we mentioned before, we think that this is at least as challenging as showing convergence of $\Delta_t$.
To our knowledge, there is no easy fix to get a result as the one claimed by \citet[Lemma 4]{yen_dual_2016} or \citet[Lemma 4]{huang_greedy_2017}.

\subsection{Comparison with UniPDGrad} %
\label{sub:comparison_with_unipdgrad}

The Universal Primal-Dual Gradient Method (UniPDGrad) by~\citet{yurtsever2015universal} is a general method to optimize problem of the form,
\begin{equation}\label{eq:UniPDGrad}
  \min_{\uu \in \mathcal{C}}\{ \tilde f(\uu) \;:\: \bm{A}\uu - \bm{b} \in \mathcal{K}\} \,,
\end{equation}
where $f$ is a convex function, $\bm{A}$ is a matrix, $\bm{b}$ a vector and $\mathcal{C}$ and $\mathcal{K}$ two closed convex sets. \eqref{eq:opt} is a particular case of their framework. There exist many ways to reformulate their framework for our application, but most of them are not practical because they require too expensive oracles. If the problem,
\begin{equation}
\argmin_{\x \in \mathcal{X}} \; f(\x) + \innerProd{\y}{M\x}  
\end{equation}
is easy to compute (which is not the case in practice most of the time but happens when $f$ is linear) then we can set $\mathcal{C} = \X, \, \mathcal{K}= \{\bm{0}\}, \, \uu = \x,\, \bm{A} = M$ and $\tilde f = f$.
Otherwise, we propose the reformulation that seemed to be the most relevant, this is the reformulation used in their experiments~\citep[Eq.19 \& 41]{yurtsever2015universal}.
If we set $\mathcal{K} = \{\bm{0}\},\;\mathcal{C} = \R^p \times \X, \; \uu = (\x,\bm{r}), \; \bm{b} = \bm{0},\; \tilde f(\uu) = f(\bm{r})$ and $\bm{A}$ such that $\bm{A}\uu= (M\x,\x-\bm{r})$ we get,
\begin{equation}\label{eq:UniPDGrad2}
  \min_{\bm{r} \in \R^p,\, \x \in \mathcal{X}}\{ f(\bm{r}) \;:\: \x = \bm{r}, \; M\x = 0\} \,,
\end{equation}
which is a reformulation of~\eqref{eq:opt}.
They derive their algorithm optimizing the (negative) Lagrange dual function $g$. The Lagrange function is,
\begin{equation}
  \LL(\x,\bm{r},\y,\bm{\lambda}) := f(\bm{r}) - \innerProd{\bm{\lambda}}{\bm{r}-\x} + \innerProd{\y}{M\x}
\end{equation}
where $\bm{\lambda}$ is the dual variable associated with the constrain $\bm{r} = \x$. Then, the (negative) Lagrange dual function is, 
\begin{align}
  g (\bm{\lambda},\y) 
  &= - \min_{\x \in\X, \bm{r} \in \R^{p}}  f(\bm{r}) - \innerProd{\bm{\lambda}}{\bm{r}-\x} + \innerProd{\y}{M\x} \notag\\ 
  &= -\min_{\bm{r} \in \R^p}  f(\bm{r}) - \innerProd{\bm{\lambda}}{\bm{r}} - \min_{\x \in \X} \innerProd{M^\top \y +\bm{\lambda}}{\x} \notag \\
  & = - f^*(\bm{\lambda})  - \min_{\x \in \X} \innerProd{M^\top \y +\bm{\lambda}}{\x} \,.
\end{align}
Their algorithm optimizes this dual function.
Computing the subgradients of the function $g$ requires to compute the Fenchel conjugate of $f$ and a LMO.

Note that \FWAL\ does not require the efficient computation of the Fenchel conjugate.

UniPDGrad computes different updates than \FWAL\ and require different assumptions for the theoretical guaranties. Particularly, \citet{yurtsever2015universal} assume the Hölder continuity of the dual function $g$. Since, in practice, the \LMO\ is not better than $0$-Hölder continuous (i.e. has bounded subgradient), we have to also assume that $f^*$ has bounded subgradients to insure the $0$-Holder continuity of the dual function. By duality, if the subgradients of $f^*$ are bounded then the support of $f$ is bounded. It is a strong assumption if we want to be able to compute the Fenchel conjugate of $f$.
Nevertheless, it seems that their proof could be extended to a dual function $g$ written as a sum of Hölder continuous functions with different Hölder continuity parameters. It would extend UniPDGrad convergence result to $f$ strongly convex ($f^*$ $1$-Hölder continuous).

In terms of rate both algorithms are hard to compare since the assumptions are different but in any case the analysis of UniPDGrad does not provide a geometric convergence rate when the constraint set $\X$ is a polytope (and $f$ a generalized strongly convex function).

It remains an open question to explore more in details and compare all the possible reformulation of~\eqref{eq:UniPDGrad} to optimize~\eqref{eq:opt} with UniPDGrad.
\section{Technical results on the Augmented Lagrangian formulation} %
\label{sec:technical_lemmas_on_L_and_d}
Let us recall that the Augmented Lagrangian function is defined as
\begin{equation}\label{eq:recall_L}
  \LL(\x,\y) :=  f(\x) +\mathbf{1}_{\X}(\x) + \innerProd{\y}{M\x} + \tfrac{\lambda}{2} \|M\x\|^2\;, \quad \forall (\x,\y) \in \R^m \times \R^d\;,
\end{equation}
where $ f$ is an $L$-smooth function, $\mathbf{1}_{\X}$ is the indicator function over the convex compact set $\X := \X_1 \times \ldots \times \X_K \subset \R^m$, $M$ is the matrix defined in~\eqref{eq:opt_reformulation}, and $m = d_1 + \ldots + d_K$. The augmented dual function $d$ is 
\begin{equation}\label{eq:def_dual_func}
d(\y) := \min_{\x \in \X} \LL(\x,\y)\,.
\end{equation}
Strong duality ensures that $\X^* \times \Y^*$ is the set of saddle points of $\LL$ where $\X^*$ is the optimal set of the primal function $p$ defined as, 
\begin{equation}\label{eq:def_primal_func}
p(\x) := \max_{\y \in \R^d} \LL(\x,\y)
\end{equation}
and $\Y^*$ is the optimal set of $d$.
 In this section we will first prove that the augmented dual function is smooth and have a property similar to strong convexity around its optimal set. It will be useful for subsequent analyses to detail the  properties of the augmented Lagrangian function $\LL$.

\subsection{Proof of Theorem~\ref{thm:dual_function_PL}} %
\label{sub:properties_of_the_dual_function_}
In this section we prove Theorem~\ref{thm:dual_function_PL}. We start with some properties of the dual function $d$.
This function can be written as the composition of a linear transformation and the Fenchel conjugate of 
\begin{equation}
f_\lambda(\x) :=  f(\x) + \frac{\lambda}{2}\|M\x\|^2 + \mathbf{1}_\X(\x) \,,
\end{equation}
where $\mathbf{1}_{\X}$ is the indicator function of $\X$. More precisely, if we denote by $^\star: f \mapsto f^*$ the Fenchel conjugate operator, then we have,
\begin{equation}\label{eq:recall_d}
  d(\y)
  := \min_{\x \in \R^m} \LL(\x,\y)
  = - \max_{\x \in \R^m} \innerProd{-M^\top \y}{\x} - f_\lambda(\x)
  = - f_\lambda^\star(-M^\top \y) \;.
\end{equation}

\paragraph{Smoothness of the augmented dual function.} %
\label{par:smoothness_of_the_augmented_dual_function_}

The smoothness of the augmented dual function is due to the duality between strong convexity and strong smoothness~\citep{rockafellar1998variational}. In order to be self-contained, we provide the proof of this property given by \citet{hong_linear_2017}.

\begin{proposition}[Lemma 2.2~\citep{hong_linear_2017}]\label{prop:dual_smooth}
If $ f$ is convex, the dual function $d$~\eqref{eq:def_dual_func} is $1/\lambda$-smooth, i.e.,
\begin{equation}
  \nabla d(\y) = M \hat \x(\y), \quad \text{where} \quad \hat \x(\y) \in \argmin_{\x \in \X} \LL(\x,\y) \;, \; \forall \y \in \R^d\;,
\end{equation}
and 
\begin{equation}
  \|\nabla d(\y) - \nabla d(\y') \| \leq \frac{1}{\lambda}\|\y- \y'\| \,\quad \forall \,\y,\y' \in \R^d\,.
\end{equation}
\end{proposition}
\proof
We will start by showing that the quantity $M \hat \x(\y)$ has the same value for all $\hat\x(\y) \in \argmin_{\x \in \X} \LL(\x,\y)$. We reason by contradiction and assume there exists $\x, \x' \in \argmin_{\x \in \X} \LL(\x,\y)$ such that $M\x \neq M\x'$. Then by convexity of $f$ and strong convexity of $\|\cdot\|^2$ we have that
\begin{equation}
 d(\y) =  \frac{1}{2} \LL(\x, \y) + \frac{1}{2} \LL(\x',\y) >  f(\bar \x) + \innerProd{\y}{M\bar \x} + \frac{\lambda}{2} \|M\bar \x\|^2 = \LL(\bar \x,\y) \;,
\end{equation}
where $\bar \x := \frac{\x + \x'}{2}$ and the inequality is strict because we assumed $M \x \neq M\x'$. This contradict the assumption that $\x, \x' \in \argmin_{\x \in \X} \LL(\x,\y)$. 
To conclude, \citet{danskin_directional_1967}'s Theorem claims that $\partial d(\y) = \{  M \hat \x(\y), \, | \, \hat \x(\y) \in \argmin_{\x \in \X} \LL(\x,\y) \}$ which is a singleton in that case. 
The function $d$ is then differentiable.

For the second part of the proof, let $\y,\y' \in \R^m$ and let $\x, \x' \in \X$ be two respective minimizers of $\LL(\cdot,\y)$ and $\LL(\cdot,\y')$. Then by the first order optimality conditions we have
\begin{equation}
  \innerProdCompressed{\nabla f(\x) + M^\top \y + \lambda M^\top M \x}{\x' - \x} \geq 0,
  \; \;
   \innerProdCompressed{\nabla f(\x') + M^\top \y' + \lambda M^\top M \x'}{\x - \x'} \geq 0 \;.
\end{equation}
Adding these two equation gives,
\begin{equation}
  \innerProdCompressed{\nabla f(\x) - \nabla f(\x') + M^\top (\y - \y') + \lambda M^\top M (\x - \x')}{\x' - \x} \geq 0 \;,
\end{equation}
but since $f$ is convex, $\innerProd{\nabla f(\x) - \nabla f(\x')}{\x -\x'} \geq 0$, and so
\begin{equation}
   \innerProd{   \y - \y'}{M(\x' - \x)} \geq  -\lambda  \innerProd{  M (\x - \x')}{M(\x' - \x)} \; .
\end{equation}
Finally, by the Cauchy-Schwarz inequality, we have
\begin{equation}
  \|\y - \y' \| \geq \lambda \|M\x - M\x'\| = \lambda  \|\nabla d(\y) - \nabla d(\y') \| \, .
\end{equation}
\endproof

\paragraph{Error bound on the augmented dual function.} %
\label{par:pl_property_of_}
After having proved that the dual function is smooth, we will derive an error bound~\citep{pang_error_1997,pang_posteriori_1987} on this function. Error bounds are related the Polyak-\L{}ojasiewic (PL) condition first introduced by~\citet{polyak_gradient_1963} and the same year in a more general setting by~\citet{lojasiewicz_topo_1963}. Recently, convergence under this condition has been studied with a machine learning perspective by~\citet{karimi_linear_2016}.

Recall that, in this section, our goal is to prove Thm.~\ref{thm:dual_function_PL}. We start our proof with lemma using the smoothness of $\LL$. 
\begin{lemma} Let $d$ be the augmented dual function~\eqref{eq:recall_d}, if $f$ is a $L$-smooth convex function and $\X$ a compact convex set, then for all $\y \in \R^d$ and $\y^* \in \Y^*$,
\begin{equation}
     d^* - d(\y)
    \geq \frac{1}{2L_\lambda D^2}\min \left( \max_{\x\in\X}\innerProd{\y^*-\y}{M\x}^2,L_\lambda D^2\max_{\x\in\X}\innerProd{\y^*-\y}{M\x}\right)
\end{equation}
where $D := \max_{(\x,\x')= \in \X^2}$ is the diameter of $\X$ and $L_\lambda := L + \lambda \|M^\top M\| $.
\end{lemma}
\proof
Let us consider $\x \in \X \subset \R^p$, $\bm{n} \in \partial f_\lambda(\x)$ a subgradient of $f_\lambda$ and the function $g_{\x}$ defined as:
\begin{equation}\label{eq:def_g}
  g_{\x}(\uu) := f_\lambda(\uu + \x) - f_\lambda(\x) - \innerProd{\uu}{\bm{n}} \;, \forall \uu \in \R^p\;.
\end{equation}
Since $f + \frac{\lambda}{2} \|M\cdot\|^2$ is $L_\lambda$-smooth, we have that $g_{\x}(\uu) \leq \frac{L_\lambda}{2} \|\uu\|^2 + \mathbf{1}_{\X}(\uu + \x) =: h_{\x}(\uu), \; \forall \uu \in \R^m$.  
By standard property of Fenchel dual (see for instance~\citep[Lemma 19]{shalev-shwartz_equivalence_2010}) we know that 
\begin{equation} \label{eq:proof_g_star_leq_h_star}
g_{\x}(\uu) \leq h_{\x}(\uu),\; \forall \uu \in \R^m \Rightarrow g_{\x}^\star(\vv) \geq h^\star_{\x}(\vv)\, , \; \forall  \vv \in \mathcal \R^m \,. 
\end{equation}
Dual computations give us for all $\vv$,
  \begin{align}
    g_{\x}^\star(\vv)
    & = \max_{\uu\in \R^m} \left[ \innerProd{\uu}{\vv} - f_\lambda(\uu + \x) + \innerProd{\uu}{\bm{n}} \right]+  f_\lambda(\x) \notag \\
    & = \max_{\uu\in \R^m} \left[ \innerProd{\uu}{\vv + \bm{n}} - f_\lambda(\uu + \x) \right]+  f_\lambda(\x)  \notag \\
    & = f^\star_\lambda( \vv + \bm{n}) + f_\lambda(\x) - \innerProd{\x}{\vv  + \bm{n}} \notag\\
    & = f^\star_\lambda(\vv  + \bm{n}) - f^\star_\lambda(\bm{n}) - \innerProd{\x}{\vv} \;,\label{eq:deriv_g_star}
  \end{align}
where in he last line we used that $\forall \bm{n} \in \partial f_\lambda(\x), \; \innerProd{\x}{\bm{n}} = f_\lambda(\x) + f_\lambda^\star(\bm{n})$ (for a proof see for instance, \citep[Lemma 17]{shalev-shwartz_equivalence_2010}).

By strong duality we have that $\X^* \times \Y^*$ is the set of saddle points, where $\X^*$ and $\Y^*$ are respectively the optimal sets of $p(\cdot)$ and $d(\cdot)$, respectively introduced in \eqref{eq:def_primal_func} and \eqref{eq:def_dual_func}.
In the following we will fix a pair $(\x^*, \y^*) \in \X^*\times \Y^*$.
Then by the stationary conditions we have
\begin{equation}\label{eq:stationary_conditions}
    M^\top \y^* \in - \partial f_\lambda(\x^*), \quad \text{and} \quad M\x^* = 0\;.
\end{equation}
Equivalently, there exist $\bm{n} \in \partial f_\lambda(\x^*)$ such that
\begin{equation}
\label{eq:nabla_f_x_star_cone_M_y_star}
\bm{n} =- M^\top \y^* \;. %
\end{equation}
For all $\y^* \in  \Y^*$ we can set $\x = \x^*$ and $\bm{n} \in \partial f_\lambda(\x^*)$ such that $\bm{n} =- M^\top \y^*$ in \eqref{eq:def_g} to get the following inequality,
\begin{align}
  d^* - d(\bm v + \y^*)
    &\; = \; f^\star_\lambda(-M^\top  \vv -M^\top\y^*) - f^\star_\lambda(-M^\top\y^*) \notag \\
    &\stackrel{\eqref{eq:nabla_f_x_star_cone_M_y_star}}{=} f^\star_\lambda(-M^\top\vv +  \bm{n}) - f^\star_\lambda(\bm{n}) \notag \\
    &\stackrel{\eqref{eq:stationary_conditions}}{=} f^\star_\lambda(-M^\top\vv +  \bm{n}) - f^\star_\lambda(\bm{n}) - \innerProd{\x^*}{-M^\top\vv} \notag \\
    &\stackrel{\eqref{eq:deriv_g_star}}{=} g^\star(-M^\top\vv) \notag \\
    &\stackrel{\eqref{eq:proof_g_star_leq_h_star}}{\geq} \; h^\star_{\x^*}(-M^\top  \vv) \;, \quad \forall \vv \in \R^d\;, \label{eq:proff_dual_h_star}
\end{align}
where for all $\vv \in \R^d$,
\begin{align}
  h^\star_{\x^*}(- M^\top \vv)
  &:= \max_{\x \in \R^m}[\innerProd{\x}{-M^\top\vv} - h_{\x^*}(\x)]\\
  & = \max_{\x \in \R^m}[\innerProd{\x}{-M^\top\vv} -  \frac{L_\lambda}{2} \|\x\|^2 - \mathbf{1}_{\X}(\x + \x^*)] \\
  & = \max_{\x+\x^* \in \X}[\innerProd{\x}{-M^\top\vv} -  \frac{L_\lambda}{2} \|\x\|^2] \label{eq:computation_h_star}
\end{align}
Let us choose $\y \in \R^d$ and set $\vv = \y - \y^*$, where $\y^* = P_{\Y^*}(\y)$. Then combining \eqref{eq:proff_dual_h_star} and \eqref{eq:computation_h_star} we get for all $\x\in \X$, and $\gamma \in [0,1]$ that $\gamma \x +(1- \gamma) \x^* \in \X$ and then,
\begin{align}
  d^* - d(\y)
    & \geq  \left[-\gamma \innerProd{M^\top(\y-\y^*)}{\x -\x^*} - \frac{L_\lambda}{2} \gamma^2 \|\x-\x^*\|^2 \right] \\
    & \geq \frac{1}{2} \left[2\gamma \innerProd{\y-\y^*}{-M\x} - \gamma^2 L_\lambda D^2 \right]  \label{eq:d_gamma}\;,
\end{align}
where $D := \max_{(\x,\x') \in \X^2}\|\x-\x'\|$ is the diameter of $\X$. Since $d^* \geq d(\y)$ the last equation can give a non trivial lower bound when $\max_{\x\in\X} \innerProd{\y-\y^*}{-M\x} >0$, we will now prove that is it always the case when $\y \notin \Y^*$.

In this proof, for $\x \in \X$ we note $N_c^{\X}(\x)$ the normal cone to $\X$ at $\x$ defined as 
\begin{equation}
  N_c^{\X}(\x) := \{ \uu \in \R^m \; |\; \innerProd{\uu}{\x -\x'} \geq 0\, , \; \forall \x' \in \X \}
\end{equation}
the reader can refers to~\citep{bauschke2011convex} for more properties on the normal cone.
If $\y \notin \Y^*$, then the necessary and sufficient stationary conditions lead to (recall that $M \x^* = 0$)
\begin{equation}
  \nabla f(\x^*) + M^\top \y \notin -N_c^{\X}(\x^*)\;,
\end{equation}
that is, there exist $\x \in \X$ such that $\innerProd{\nabla f(\x^*) + M^\top\y}{\x-\x^*} < 0$. Using \eqref{eq:nabla_f_x_star_cone_M_y_star} gives
 \begin{equation}
  \begin{aligned}
    0
    & \;>\; \innerProd{\nabla f(\x^*) + M^\top\y}{\x-\x^*} \\
    & \stackrel{\eqref{eq:nabla_f_x_star_cone_M_y_star}}{=} \innerProd{-M^\top \y^* - \uu + M^\top\y}{\x-\x^*} \\
    & \;\geq\;  \innerProd{\y-\y^*}{M\x} \;,
  \end{aligned}
 \end{equation}
where for the last inequality we use the fact that $\uu \in N_c^{\X}(\x^*)$ and $M\x^* = 0$.
Then we have
\begin{equation}\label{eq:gap_positive}
  \max_{\x\in\X} \innerProd{\y-\y^*}{-M\x} >0, \quad \forall \y \notin \Y^*\;.
 \end{equation}
Optimizing Eq.~\eqref{eq:d_gamma} with respect to $\gamma \in [0,1]$ we get the following:
\begin{itemize}
  \item If $0<\max_{\x\in\X}\innerProd{\y^*-\y}{M\x} \leq L_\lambda D^2 $, the optimum of \eqref{eq:d_gamma} is achieved for $\gamma= \frac{\max_{\x\in\X}\innerProd{\y^*-\y}{M\x}}{ L_\lambda D^2} \leq 1$ and we have,
   \begin{equation}
  d^* - d(\y) \geq \frac{1}{2L_\lambda D^2} \max_{\x\in\X}\innerProd{\y^*-\y}{M\x }^2  \;,
\end{equation}
  \item Otherwise, if $\max_{\x\in\X}\innerProd{\y^*-\y}{M\x}> L_\lambda D^2$, the optimum of \eqref{eq:d_gamma} is achieved for $\gamma =1$, giving
  \begin{equation}\label{eq:lb_dual_directions_1}
  d^* - d(\y)
    \geq \frac{1}{2}\max_{\x\in\X}\left[2 \innerProd{\y^*-\y}{M\x} - L_\lambda D^2 \right]
    \geq \frac{1}{2}\max_{\x\in\X}\innerProd{\y^*-\y}{M\x}\;.
\end{equation}
\end{itemize}
Combining both cases leads to
\begin{equation}\label{eq:dual_subopt_gap}
     d^* - d(\y)
    \geq \frac{1}{2L_\lambda D^2}\min \left( \max_{\x\in\X}\innerProd{\y^*-\y}{M\x}^2,L_\lambda D^2\max_{\x\in\X}\innerProd{\y^*-\y}{M\x}\right) \,.
\end{equation}
\endproof
Since our goal is to get an error bound on the dual function $d$ we divide and multiply by $\|\y-\y^*\|$ the quantities $\max_{\x\in\X}\innerProd{\y^*-\y}{M\x}$ in \eqref{eq:dual_subopt_gap}, making appear the desired norm and a constant $\alpha$ defined as
\begin{equation}\label{eq:definition_alpha}
  \alpha := \inf_{\substack{\y \in \R^d\setminus\Y^*\\
          \y^*= P_{\Y^*}(\y)}} \sup_{\x\in\X}\innerProd{\frac{\y^*-\y}{\|\y^*-\y\|}}{M\x} \,.
\end{equation}
Recall that $\y^* := P_{\Y^*}(\y)$ and consequently $\|\y - \y^*\| = \dist(\y,\Y^*)$.
Our goal is now to show that $\alpha>0$.
\paragraph{Proof that $\alpha$ is positive.} %
\label{par:proof_that_}
In order to prove that $\alpha$ is positive we need to get results on the structure of $\Y^*$. 
First, let us start with a topological lemma,
\begin{lemma}\label{lemma:0_in_int}
Let $(C_k)_{k\in[K]}$ be a collection of nonempty convex sets. We have that 
\begin{equation}\label{eq:0_in_int}
  0 \in \relint(C_k) \, , \, k \in [K] \Rightarrow 0 \in \relint\left(\underset{k=1}{\overset{K}{+}} C_k\right) \,.
\end{equation}
\end{lemma}
\proof
In order to prove this result we will prove two intermediate results. Recall that the cone $\cone(C)$ generated by a convex set $C$ is defined as 
\begin{equation}
  \cone(C) := \{ \lambda \x \;:\; \x \in C\}\,.
\end{equation}
For more details on the topological properties of the convex set set for instance~\citep{rockafellar_convex_1970}.
\begin{itemize}
  \item 
  The first one is a characterization: 
  \begin{equation} \label{eq:0_relint_cone_equal_span}
    0 \in \relint(C) \Leftrightarrow \cone(C) = \Span(C) \,.
  \end{equation}
  $\Rightarrow$: Let $\x \in C$,
  \begin{align*}
  \x \in \Span(C) 
  &\Rightarrow \exists \, \lambda_i \in \R ,\, \x_i \in C ,\, i \in \{1,\ldots,n\} \quad \text{s.t.} \quad \x = \sum_{i=1}^n \lambda_i \x_i  \\
  &\Rightarrow \exists \, \lambda_i \in \R ,\, \x_i \in C ,\, i \in \{1,\ldots,n\} \quad \text{s.t.} \quad \x = \lambda \sum_{i=1}^n \frac{\lambda_i \x_i}{\lambda} \,,\; \lambda >0 \\
  & \Rightarrow \exists \lambda >0 \,, \;  \tilde \x_i \in C\, ,\; i \in \{1,\ldots,n\} \quad \text{s.t.} \quad  \x = \lambda \sum_{i=1}^n \tilde \x_i  \\
  & \Rightarrow \x \in \cone(C)\,.
  \end{align*}
  where the last line is due to the fact that for $\lambda$ small enough $\frac{\lambda_i \x_i}{\lambda} \in C$ because $ 0 \in \relint(C)$.

  By definition we have that $\cone(C) \subset \Span(C)$. 

  Then, we have proved that $ 0 \in \relint(C) \Rightarrow \cone(C) = \Span(C) $

  $\Leftarrow$:
  If $\relint(C) = \{0\}$, then $\{0\} = \relint(C) = \cone(C) = \Span(C)$.

  Otherwise, let $\x \in \relint(C)\setminus\{0\}$, using our hypothesis we have that, 
  \begin{equation}
  -\x \in \Span(C) = \cone(C) \Leftrightarrow 0 \in \cone(C) + \x \,.   
  \end{equation}
  Then there exist $\x' \in C$ and $\lambda >0$ such that, 
  \begin{equation}
  0 = \lambda \x' + \x \Leftrightarrow 0 =  \frac{\lambda}{1+\lambda} \x' + \frac{1}{1+\lambda} \x\,.  
  \end{equation} 
  Since $ 1\geq \frac{1}{1+\lambda} > 0$ and $\x \in \relint(C)$, we have by \cite[Theorem~6.1]{rockafellar_convex_1970} that $0 \in \relint(C)$. 

  \item The second one is a property on the sum of the convex cones generated by $(C_k)$:
  \begin{equation} \label{eq:sum_cone_cone of sum}
     0 \in \relint(C_k)\,,\; k \in \{1,\ldots,K\} \Rightarrow \underset{k=1}{\overset{K}{+}}\cone(C_k) = \cone\left(\underset{k=1}{\overset{K}{+}}C_k\right) \,.
  \end{equation}
  Let, then,
  \begin{align*}
  \x \in \underset{k=1}{\overset{K}{+}}\cone(C_k)
  & \Leftrightarrow
  \exists \tilde \x_k \in \cone(C_k), \; k \in \{1,\ldots,K\} \quad \text{s.t.} \quad \x = \sum_{k=1}^k  \tilde \x_k \\
  & \Leftrightarrow
  \exists \x_k \in C_k, \; \lambda_k \in \R \,,\; k \in \{1,\ldots,K\} \quad \text{s.t.} \quad \x = \sum_{k=1}^k \lambda_k \x_k \\
  & \Leftrightarrow  
  \exists \x_k \in C_k, \; \lambda_k \in \R \,,\; k \in \{1,\ldots,K\} \quad \text{s.t.} \quad \x = \lambda\sum_{k=1}^k \frac{\lambda_k \x_k}{\lambda} \,,\; \lambda>0 \\
  & \Leftrightarrow
  \x \in \cone\left(\underset{k=1}{\overset{K}{+}}C_k\right) \,.
  \end{align*}
  For the last equivalence we used that $0 \in \relint(C_k)\,,\; k \in \{1,\ldots,K\}$.
\end{itemize}
Now we can prove our lemma using~\eqref{eq:0_relint_cone_equal_span} and~\eqref{eq:sum_cone_cone of sum}:
\begin{align*}
   0 \in \relint(C_k) \; , \; k \in [K]
   & \stackrel{\eqref{eq:0_relint_cone_equal_span}}{\Rightarrow}  \cone(C_k) = \Span(C_k) \; , \; k \in [K] \\
   0 \in \relint(C_k) \; , \; k \in [K]
   & \stackrel{\eqref{eq:sum_cone_cone of sum}}{\Rightarrow}
   \cone\left(\underset{k=1}{\overset{K}{+}}C_k\right) =  \underset{k=1}{\overset{K}{+}}\cone(C_k) = \underset{k=1}{\overset{K}{+}}\Span(C_k)  = \Span\left(\underset{k=1}{\overset{K}{+}}C_k\right) \\
   & \stackrel{\eqref{eq:0_relint_cone_equal_span}}{\Rightarrow} 
   0 \in \relint\left(\underset{k=1}{\overset{K}{+}} C_k\right) 
\end{align*}
\endproof
Let us recall the supplementary assumption needed to prove Theorem~\ref{thm:dual_function_PL}. 
\begin{repassumption}{assump:main}
$\exists \, \bar \x^{(k)} \in \relint(\X_k) ,\, k \in \{1,\ldots,K\},\, s.t.,\,\sum_{k=0}^K A_k \bar \x^{(k)} = 0$.
\end{repassumption}
This assumption is required in the proof of the following lemma, 
\begin{lemma}\label{lemma:form_Y_star}
Under Assumption~\ref{assump:main}, the optimal set $\Y^*$ of the augmented dual function $d(\cdot)$~\eqref{eq:recall_d} can be written as 
\begin{equation}
  \Y^* = \mathcal K + V\,,
\end{equation}
where $V :=   \cap_{k=1}^K \big(A_k(\Span(\X_k-\bar \x^{(k)}))\big)^\perp$ and $\mathcal K \subset V^\perp$ is a compact set.
\end{lemma}
We define $\Span(\X_k -\bar \x^{(k)})$ as the linear span of the feasible direction from $\bar \x^{(k)}$. Since $\bar \x^{(k)}$ is a relative interior point of the convex $\X_k$ we have $\Span(\X_k- \bar \x^{(k)}) = \{ \lambda(\x^{(k)} - \bar\x^{(k)}) \; : \; \x^{(k)} \in \X_k , \; \lambda>0\}$.
\proof

For any $\x^* \in \X^*$, a necessary and sufficient condition for any $\y^*$ to be in $\Y^*$ is
\begin{equation}
  \nabla f(\x^*) + M^\top \y^* \in - N_c(\x^*)\;,
\end{equation}
meaning that
\begin{equation}
  -A_k^\top \y^* \in N_c^{\X_k}(\x^*) + \nabla_{\x^{(k)}} f(\x^*)\,, \quad k \in \{1,\ldots,K\}\,.
\end{equation}
Then noting $g_k := \nabla_{\x^{(k)}} f(\x^*) + \lambda M \x$ we have the following equivalences,
\begin{align*}
    \y^* \in \Y^* 
    &\Leftrightarrow -A_k^\top \y^* \in N_c^{\X_k}(\x^*) + g_k\,, \quad k \in \{1,\ldots,K\} \\
    &\Leftrightarrow A_k^\top \y^* + g_k \in -N_c^{\X_k}(\x^*) \,, \quad k \in \{1,\ldots,K\} \\
    &\Leftrightarrow  \innerProd{-A_k^\top \y^* - g_k}{\x^{(k)} - (\x^*)^{(k)}} \leq 0 \; ;\; \forall \x^{(k)} \in \X_k \,, \quad k \in \{1,\ldots,K\} \\
    &\Leftrightarrow \innerProd{-\y^*}{A_k(\x^{(k)} - (\x^*)^{(k)})} \leq \innerProd{g_k}{\x^{(k)} - (\x^*)^{(k)}} \; ;\; \forall \x^{(k)} \in \X_k \,, \quad k \in \{1,\ldots,K\}
\end{align*}
Then we can notice that if we write $\y^* = \y_1^* + \y_2^*$ with $\y_1^* \in V:= \cap_{k=1}^{K} \big(A_k (\Span(\X_k- \bar \x^{(k)})\big)^\perp$ and $\y_2^* \in V^\perp$ we get,
\begin{equation}\label{eq:proof_Y_start_equiv}
  \y^* \in \Y^* 
    \Leftrightarrow \innerProd{-\y_2^*}{A_k(\x^{(k)} - (\x^*)^{(k)})} \leq \innerProd{g_k}{\x^{(k)} - (\x^*)^{(k)}} \; ;\; \forall \x^{(k)} \in \X_k \,, \quad k \in [K] \,.
\end{equation}
Note that there is no conditions on $\y_1^*$. 

Let us get a necessary condition on $\y_2^*$. Eq.~\eqref{eq:proof_Y_start_equiv} implies,
\begin{align*}
 \y^* \in \Y^* 
 &\Rightarrow \innerProd{-\y_2^*}{\sum_{k=1}^KA_k(\x^{(k)} - (\x^*)^{(k)})} \leq \sum_{k=1}^K \innerProd{g_k}{\x^{(k)} - (\x^*)^{(k)}} \; ;\; \forall \x^{(k)} \in \X_k \,, \quad k \in [K]   \\
 &\Rightarrow \innerProd{-\y_2^*}{\sum_{k=1}^KA_k(\x^{(k)} - (\x^*)^{(k)})} \leq \sum_{k=1}^K \|g_k\|\|\x^{(k)} - (\x^*)^{(k)}\| \; ;\; \forall \x^{(k)} \in \X_k \,, \quad k \in [K]   \\
 & \Rightarrow  \innerProd{-\y_2^*}{\sum_{k=1}^KA_k(\x^{(k)} - \bar \x^{(k)})} \leq \sum_{k=1}^K \|g_k\|\diam(X_k)\quad ;\quad \forall \x^{(k)} \in \X_k \,, \quad k \in [K] \,, \\
 & \qquad \text{($\x^*$ and $\bar \x$ are feasible, i.e., $\sum_{k=1}^KA_k\bar\x^{(k)} = \sum_{k=1}^KA_k (\x^*)^{(k)} =0$)} 
\end{align*}
where $\bar \x^{(k)} \in \relint(\X_k) ,\; k \in [K]$ and $M\bar \x = 0$ (Assump.~\ref{assump:main}). 
Moreover, since $V:= \cap_{k=1}^{K} \big(A_k (\Span(\X_k- \bar \x^{(k)}))\big)^\perp$ we have that $V^\perp = \overset{K}{\underset{k=1}{+}} A_k (\Span(\X_k- \bar \x^{(k)})).$ 
Then by Lemma~\ref{lemma:0_in_int}, 
\begin{align*}
\bar \x \in \relint(\X)
&\Rightarrow 0 \in \relint(A_k (\X_k- \bar \x^{(k)})) \quad k \in \{1,\ldots,K\} \\
&\stackrel{\eqref{eq:0_in_int}}{\Rightarrow} 0 \in \relint\left(\overset{K}{\underset{k=1}{+}} A_k (\X_k- \bar \x^{(k)})\right)\,,
\end{align*}
and consequently, there exists $\delta >0$ such that for all $\y_2^* \in \overset{K}{\underset{k=1}{+}} A_k (\Span(\X_k- \bar \x^{(k)})) $, we can set $\x^{(k)} \in \X_k$ such that $\sum_{k=1}^KA_k(\x^{(k)} - \bar \x^{(k)}) = -\delta \y_{2}^*/\|\y_2^*\|$. Finally, we get that,
\begin{equation}
   \y^* \in \Y^* 
 \Rightarrow  \delta\innerProd{\y_2^*}{\y_{2}^*/\|\y_2^*\|} \leq \sum_{k=1}^K \|g_k\|\diam(X_k)
 \Rightarrow \|\y_2^*\|_2 \leq \sum_{k=1}^K \frac{\|g_k\|\diam(X_k)}{\delta} \,.
\end{equation}
Thus $\mathcal{K} \subset V^\perp$ is bounded and consequently compact (because $\Y^*$ is closed).
\begin{proposition}\label{prop:description_Y_star}
If Assumption~\ref{assump:main} holds, then the set of normal directions to $\Y^*$,
\begin{equation}
 \mathcal {D} := \big\{ \bm{d} \::\: \bm{d} \in  N_c^{\Y^*}(\y^*) \text{ for $\y^* \in \Y^*$} \,,\, \|\bm{d}\| = 1 \big\}\;,
\end{equation}
is closed and consequently compact.
\end{proposition}
\proof
Let us first show that,
\begin{equation} \label{eq:normal_directions}
\mathcal {D} = \{ \y - P_{\Y^*}(\y)  \;:\; \y \in \R^d \; ; \;\|\y - P_{\Y^*}(\y)\| =1\}\,.
\end{equation}
Let $\y \in \R^d \setminus \Y^*$, by definition of the normal cone and the projection onto a convex set, we have that $\y - P_{\Y^*}(\y) \in N^{\Y^*}_c(P_{\Y^*}(\y))$. Conversely, for any $\y^* \in \Y^*$ and $\dd \in  N_c^{\Y^*}(\y^*)$ such that $\|\dd\|=1$, we have that $\y^* = P_{\Y^*}(\y^*+\dd)$ and $\y^* + \dd \notin \Y^*$.

With the same notation as Lemma~\ref{lemma:form_Y_star}, we can write $\y \in \R^d$ a unique way as $\y = \y_1+\y_2$ where $\y_1 \in V$ and $\y_2 \in V^\perp$. Then since $\Y^* = V \overset{\perp}{+} \mathcal{K}$ we get that $P_{\Y^*}(\y) = \y_1 + \bm{\kappa}$ where $\bm{\kappa}\in \mathcal{K}$. Then $\y - P_{\Y^*}(\y) = \y_2 - \bm{\kappa}$ where $P_{\mathcal{K}}(\y_2)=  \bm{\kappa}$. 
Conversely, for any couple $(\y_2,\kappa) \in V^\perp \times \mathcal K$ such that $P_{\mathcal{K}}(\y_2)=  \bm{\kappa}$, we have that $\y_2 - \bm{\kappa} \in N_c^{\mathcal{K}}(\bm{\kappa})$.

If we call $\phi : \y \mapsto \y-P_{\mathcal{K}}(\bm{\y})$, then $\mathcal{D} = \phi(A) $
where $A = \{\y_2 \in V^\perp \; ;\; \dist(\y_2,\mathcal{K}) =1\}$ is a compact (because $\mathcal{K}$ is compact). Then since $\phi$ is continuous, $\mathcal{D}$ is a compact.
\endproof

Now we can apply this result to bound the $\alpha$ constant introduced in Eq.~\eqref{eq:definition_alpha}. We notice that using~\eqref{eq:normal_directions}, we can write that definition as
\begin{equation}\label{eq:re_def_alpha}
\alpha = \inf_{\substack{\y \in \R^d\setminus\Y^*\\
            \y^*= P_{\Y^*}(\y)\\
            \bm{d} = \y^*\! - \y,\,
            \|\bm{d}\|=1}
        } \sup_{\x\in\X}\innerProd{\bm{d}}{M\x} \;.
\end{equation}
The function $ \bm{d} \mapsto \sup_{\x\in\X}\innerProd{\bm{d}}{M\x}$ is convex
(as a supremum of convex function) and then is continuous on the interior of its domain which is $\R^d$ because $\X$ is bounded.
Since $\mathcal{D}$ is compact, the infimum is achieved. Then, there exist $\y \in (\R^d\setminus \Y^*)\times \Y^*$ such that, $\y^* = P_{\Y^*}(\y)$, $\|\y^*-\y\|=1$ and, 
\begin{equation}
\alpha = \max_{\x\in \X}\innerProd{\y^* - \y}{M\x}.  
\end{equation}
By Equation~\eqref{eq:gap_positive}, since $\y$ is non optimal, we conclude that $\alpha>0$.

\paragraph{Proof of Thm.~\ref{thm:dual_function_PL} and a Corollary.} %
\label{par:proof_of_thm_thm:dual_function_pl_and_a_corollary_}

\begin{reptheorem}{thm:dual_function_PL}
 Let $d$ be the augmented dual function~\eqref{eq:recall_d}, if $f$ is a $L$-smooth convex function and $\X$ a compact convex set and if Assumption~\ref{assump:main} holds, then for all $\y \in \R^d$ there exist a constant $\alpha>0$ such that,
\begin{equation}\label{eq:lb_dual_directions_2}
  d^* - d(\y)
  \geq \frac{1}{2L_\lambda D^2}\min \left\{
        \alpha^2\dist(\y,\Y^*)^2, \alpha L_\lambda D^2\dist(\y,\Y^*)
  \right\} \;,
\end{equation}
where $D := \max_{\x,\x' \in \X}\|\x-\x'\|$ is the diameter of $\X$.
 \end{reptheorem}
 \proof Recall that we proved 
 \begin{equation}
     d^* - d(\y)
    \geq \frac{1}{2L_\lambda D^2}\min \left( \max_{\x\in\X}\innerProd{\y^*-\y}{M\x}^2,L_\lambda D^2\max_{\x\in\X}\innerProd{\y^*-\y}{M\x}\right)~,
\end{equation}
and that $\alpha$ defined in \eqref{eq:definition_alpha} was positive~\eqref{eq:alpha_positive}. Then for all $\y \notin \Y^*$,
\begin{equation}
    d^* - d(\y)
    \geq \frac{1}{2L_\lambda D^2}\min \left( \alpha^2\dist(\y,\Y^*)^2,L_\lambda D^2\alpha\dist(\y,\Y^*)\right) \;.
\end{equation}
The same result is trivially true for $\y \in \Y^*$ (since in that case we have $d(\y) = d^*$).
 \endproof
 This Theorem leads to an immediate corollary on the norm of the gradient of $d$.
 \begin{corollary}\label{cor:lower_bound_gradient}
Under the same assumption as Theorem~\ref{thm:dual_function_PL}, for all $\y \in \R^d$ there exist a constant $\alpha$ such that,
   \begin{equation}\label{eq:lb_dual_directions}
  \|\nabla d(\y)\| \geq \frac{1}{2L_\lambda D^2}\min\{ \alpha^2\dist(\y,\Y^*),\alpha L_\lambda D^2  \}
  \;\;\;\text{and} \;\;\;
  \|\nabla d(\y)\| \geq \frac{\alpha}{\sqrt{2L_\lambda D^2}}\min\Big\{\sqrt{ d^*-d(\y)}, \sqrt{\tfrac{L_\lambda D^2}{2} }  \Big\}\,.
  \end{equation}
 \end{corollary}
 \proof
 We just need to notice that by concavity of $d$ for all $\y^* \in \Y^*$, the suboptimality is upper bounded by the linearization of the function: 
 \begin{equation} \label{eq:proof_bound_subopt}
 d^* - d(\y) \leq \innerProd{\y^* - \y}{\nabla d(\y)} \leq \dist(\y,\Y^*)\|\nabla d(\y)\| \,.   
 \end{equation}
 Then combining it with Theorem~\ref{thm:dual_function_PL} we get, 
 \begin{equation}\label{eq:proof_bound_nabla}
   \frac{1}{2L_\lambda D^2}\min \left\{
        \alpha^2\dist(\y,\Y^*), \alpha L_\lambda D^2
  \right\} \leq \|\nabla d(\y)\| \,.
 \end{equation}
This equation is equivalent to 
\begin{equation}\label{eq:proof_bound_nabla_2}
  \frac{1}{2L_\lambda D^2}
        \alpha^2\dist(\y,\Y^*) \leq \|\nabla d(\y)\|
  \quad \text{or} \quad \frac{\alpha}{2} \leq \|\nabla d(\y)\| \,.
\end{equation}
Combining the first inequality of~\eqref{eq:proof_bound_nabla_2} with \eqref{eq:proof_bound_subopt} we get,
\begin{equation}
  d^* - d(\y) \leq \frac{2 L_\lambda D^2}{\alpha^2}  \|\nabla d(\y)\|^2
  \quad \text{or} \quad \frac{\alpha}{2} \leq \|\nabla d(\y)\| \,,
\end{equation}
which is equivalent to  
\begin{equation}
  \|\nabla d(\y)\| 
  \geq  \frac{\alpha}{\sqrt{2L_\lambda D^2}}\min\{\sqrt{ d^*-d(\y)}, \sqrt{{L_\lambda D^2}/{2}}  \}\; .
\end{equation}
 \endproof

\subsection{Properties of the function $\LL$} %
\label{sub:lemmas_about_}
We will first prove that for any $\y \in \R^d$, the function $\LL(\cdot,\y)$ has a property similar to strong convexity respect to the variable $M\x$: if $\LL(\x,\y)$ is close to its minimum with respect to $\x$, then  $M\x$ is close to the image by $M$ of the minimizer of $\LL(\cdot,\y)$. More precisely,
\begin{proposition}\label{prop:strong_conv_Mx} for all $\x \in \X$ and $\y \in \R^d$, if $f$ is convex,
\begin{equation} \label{eq:strong_conv_Mx}
\|M\x-M \hat\x(\y)\|^2 \leq  \frac{2}{\lambda} \left( \LL(\x,\y) - \LL(\hat \x(\y),\y) \right)
\quad \text{where} \quad
\hat\x(\y) \in \argmin_{\x \in \X} \LL(\x,\y) \;,
\end{equation}
and $\LL(\x,\y) := f(\x) +\mathbf{1}_{\X}(\x) + \innerProd{\y}{M\x} + \tfrac{\lambda}{2} \|M\x\|^2$
\end{proposition}
\proof
By convexity of $f$ we have that,
\begin{equation}
  f(\x) - f(\hat \x(\y)) \geq \innerProd{\nabla f(\hat \x(\y))}{\x - \hat \x(\y)}\;,
\end{equation}
then by simple algebra (noting $\hat  \x = \hat \x(\y)$),
\begin{align}
   \LL(\x,\y) - \LL(\hat \x,\y)
   &\geq \innerProd{\nabla f (\hat \x)+ M^\top \y + \lambda M^\top M \hat\x}{\x - \hat \x} + \frac{\lambda}{2} \|M\x - M \hat \x\|^2 \\
   & = \innerProd{\nabla_{\x} \LL (\hat \x,\y)}{\x - \hat \x} + \frac{\lambda}{2} \|M\x - M \hat \x\|^2 \\
   &\geq \frac{\lambda}{2} \|M\x - M \hat \x\|^2 \;.
\end{align}
The last inequality come from the first order optimality condition on $\LL(\cdot,\y)$.
\endproof
Now let us introduce the key property allowing us to insure that $\x_{t}$ actually converge to $\x^*$. This proposition states that the primal gap $\Delta_t^{(p)}$ upper-bounds the squared distance to the optimum.
\begin{proposition} \label{prop:delta_bounds_input} If $f$ is a $\mu$-strongly convex function then, $\X^* = \{\x^*\}$ and we have for all $t\geq 0$,
\begin{equation}
  \frac{\mu}{2}\|\x_{t+1}-\x^*\|^2 \leq \Delta_t^{(p)}  - \Delta_t^{(d)} + \max \left( 2\Delta_t^{(d)}, \sqrt{2L_\lambda D^2 \Delta_t^{(d)}} \right)\;.
\end{equation}
and also
\begin{equation}
   \frac{\mu}{2}\|\x_{t+1}-\x^*\|^2 \leq \Delta_t^{(p)} + \frac{2L_\lambda D^2}{\alpha^2} \|M\x_{t+1}\| \|M\hat \x_{t}\|\,, \; \forall t \in \N \;;\; \dist(\y_t,\Y^*) \leq \frac{L_\lambda D^2}{\alpha}\;.
\end{equation}
\end{proposition}
\proof
We start from the identity
\begin{equation} \label{eq:subopt_eq_delta}
  f(\x_{t+1}) -  f^* = \Delta_t^{(p)} - \Delta_t^{(d)} - \innerProd{\y_{t}}{M\x_{t+1}} - \frac{\lambda}{2}\|M\x_{t+1}\|^2~.
\end{equation}
From first order optimality conditions, we get for any $\y^* \in \Y^*$ and any $\x \in \X$,
\begin{equation}
  \innerProd{\x - \x^*}{\nabla f(\x^*) + M^\top \y^* + \lambda M^\top M \x^*} \geq 0 \quad \text{and} \quad M \x^* = 0 \; ,
\end{equation}
then for $\x = \x_{t+1}$,
\begin{equation} \label{eq:first_order_cond}
  \innerProd{\y^*}{M \x_{t+1}} \geq -\innerProd{\x_{t+1} - \x^* }{\nabla f(\x^*)} \;.
\end{equation}
If $f$ is $\mu$-strongly convex, then
\begin{equation} \label{eq:strong_conv}
  -\innerProd{\x_{t+1} - \x^* }{\nabla f(\x^*)} \geq -f(\x_{t+1}) + f^* + \frac{\mu}{2}\|\x_{t+1} -\x^*\|^2 \;,
\end{equation}
then combining \eqref{eq:subopt_eq_delta}, \eqref{eq:first_order_cond} and \eqref{eq:strong_conv}, we get for any $\y^* = P_{\Y^*}(\y)$:
\begin{align}
  \frac{\mu}{2}\|\x_{t+1}-\x^*\|^2 
  &\leq \Delta_t^{(p)} - \Delta_t^{(d)} + \innerProd{\y^* - \y_{t}}{M\x_{t+1}} \\
  &\leq \Delta_t^{(p)} - \Delta_t^{(d)} + \dist(\y_t,\Y^*) \|M\x_{t+1}\|  \;. \label{eq:109}
\end{align}
Using the fact that in \eqref{eq:lb_dual_directions},
\begin{equation}
\text{either} \quad \dist(\y_t ,\Y^*) \geq \frac{L_\lambda D^2}{\alpha}
\quad \text{or} \quad
\dist(\y_t,\Y^*) \leq \frac{2L_\lambda D^2}{\alpha^2}\|\nabla d(\y_t)\| = \frac{2L_\lambda D^2}{\alpha^2}\|M \hat\x_{t}\| \;,
 \end{equation}
 leading to
\begin{equation}
  \frac{\mu}{2}\|\x_{t+1}-\x^*\|^2 \leq \Delta_t^{(p)} - \Delta_t^{(d)} + \frac{2L_\lambda D^2}{\alpha^2} \|M\x_{t+1}\| \|M\hat \x_{t}\|\,, \; \forall t \in \N \;;\; \dist(\y_t,\Y^*) \leq \frac{L_\lambda D^2}{\alpha}\;.
\end{equation}
Similarly, combining~\eqref{eq:109} and \eqref{eq:dual_subopt_gap} gives us,
\begin{equation}
  \frac{\mu}{2}\|\x_{t+1}-\x^*\|^2 \leq \Delta_t^{(p)}  - \Delta_t^{(d)} + \max \left( 2\Delta_t^{(d)}, \sqrt{2L_\lambda D^2 \Delta_t^{(d)}} \right)\;.
\end{equation}
\endproof
This property will be used to prove Theorem~\ref{thm:strong}, deducing convergence rates on $\|\x_{t+1}-\x^*\|^2$ from the convergence rates on $\Delta_t$ proved in Theorem~\ref{thm:main_general} and Theorem~\ref{thm:main_polytope}.

\section{Proof of Theorem~\ref{thm:main_general}, Theorem~\ref{thm:main_polytope} and Theorem~\ref{thm:strong}} %
\label{app:proof_of_thm}
This section is decomposed into 3 subsections. 
First, we prove some intermediate results on the sequence computed by our algorithm to get the fundamental equation~\eqref{eq:fund_descent} that we will use to prove the convergence of $(\Delta_t)_{t\in\N}$.
Then in subsection~\ref{sub:proof_of_1_in_theorem_thm:main} (respectively Subsection~\ref{sub:proof_of_theorem_thm:main_polytope}) we prove Thm.~\ref{thm:main_general} (resp. Thm.~\ref{thm:main_polytope}).
Let us recall that the Augmented Lagrangian function is defined as
\begin{equation}\label{eq:recall_L_second_time}
  \LL(\x,\y) := f(\x) +\mathbf{1}_{\X}(\x) + \innerProd{\y}{M\x} + \tfrac{\lambda}{2} \|M\x\|^2\;, \quad \forall (\x,\y) \in \R^m \times \R^d\;,
\end{equation}
where $f$ is a smooth function, $\mathbf{1}_{\X}$ is the indicator function of a convex compact set $\X \subset \R^m$. 
The augmented dual function $d$ is $d(\y) := \max_{\x \in \X} \LL(\x,\y)\,.$ The \FWAL{} algorithm computes
\begin{equation}\label{eq:app_d_update_rue}
\left\{
  \begin{aligned}
    &\x_{t+1} = \mathcal{FW}(\x_{t};\LL(\cdot,\y_{t}))\;,\\
    &\y_{t+1} = \y_{t} + \eta_{t}M \x_{t+1} \;,
  \end{aligned}
  \right.
\end{equation}
where $\mathcal{FW}(\x_{t};\LL(\cdot,\y_{t}))$ is roughly a \FW\ step from $\x_t$. (More details in App.~\ref{app:frank_wolfe_algorithms}).

\subsection{Lemma deduced from the dual variable update rule} %
\label{sub:lemma_deduced_from_the_dual_variable_update_rule}

The two following lemmas do not require any assumption on the sets or the functions, they only rely on the dual update on $\y$~\eqref{eq:app_d_update_rue}. They provide upper bounds on the decrease of the primal and the dual gaps.
They are true for all functions $f$ and constraint set $\X$.
Recall that we respectively defined the primal and the dual gap as,
\begin{equation}\label{eq:app_primal_and_dual_gaps}
  \Delta_t^{(d)}:=d^*-d(\y_{t}) \quad \text{and} \quad \Delta_t^{(p)} := \LL(\x_{t+1};\y_{t})-d(\y_{t}) \;.
\end{equation}
The first lemma upper bounds the decrease of the dual suboptimality; note that~\citet{hong_linear_2017} are probably not the firsts to provide such lemma. We are citing them because we provide the proof proposed in their paper.
\begin{lemma}[Lemma 3.2~\citep{hong_linear_2017}]\label{lemma:dual_descent}
For any $t \geq 1$, there holds
\begin{equation}
\label{eq:dual_gap}
\Delta_{t+1}^{(d)}-\Delta_{t}^{(d)}\leq -\eta_{t} \innerProd{M\x_{t+1}}{M\hat \x_{t+1}}.
\end{equation}
\end{lemma}
\proof
\begin{align}
  \Delta_{t+1}^{(d)}-\Delta_{t}^{(d)}
  & = d(\y_{t}) - d(\y_{t+1}) \notag \\
  & = \LL(\hat \x_{t}, \y_{t}) - \LL(\hat \x_{t+1}, \y_{t+1})\notag \\
  & \stackrel{(\star)}{\leq} \LL(\hat \x_{t+1}, \y_{t}) - \LL(\hat \x_{t+1}, \y_{t+1}) \label{eq:dual_gap_2}\\
  & = \innerProd{\y_{t} - \y_{t+1}}{M\hat \x_{t+1}}\notag \\
  & = -\eta_{t} \innerProd{M\x_{t+1}}{M\hat \x_{t+1}} \,, 
\end{align}
where $(\star)$ is because $\hat \x_{t}$ is the minimizer of $\LL(\cdot,\y_{t})$.
\endproof

Next we proceed to bound the decrease of the primal gap $\Delta_{t+1}^{(p)}$.
\begin{lemma}[weaker version of Lemma 3.3~\citep{hong_linear_2017}]\label{lemma:primal_descent}
Then for any $t\geq1$, we
have
\begin{equation}
\label{eq:primal_gap}
\Delta_{t+1}^{(p)}-\Delta_{t}^{(p)}\leq \eta_{t}\|M\x_{t+1}\|^2+(\LL(\x_{t+2},\y_{t+1})-\LL(\x_{t+1},\y_{t+1}))-\eta_{t} \innerProd{M\x_{t+1}}{M\hat \x_{t+1}}.
\end{equation}
\end{lemma}
\proof
We start using the definition of $\Delta_{t+1}^{(p)}$,
\begin{align*}
\Delta_{t+1}^{(p)}-\Delta_{t}^{(p)}
& \;\,=\;\, \LL(\x_{t+2},\y_{t+1}) - \LL(\hat \x_{t+1},\y_{t+1}) - \left( \LL(\x_{t+1},\y_{t}) - \LL(\hat \x_{t},\y_{t}) \right) \\
& \;\,=\;\, \LL(\x_{t+2},\y_{t+1}) - \LL( \x_{t+1},\y_{t}) + \left( \LL(\hat \x_{t},\y_{t}) - \LL(\hat \x_{t+1},\y_{t+1}) \right) \\
& \stackrel{\eqref{eq:dual_gap_2}}{\leq} \LL(\x_{t+2},\y_{t+1}) - \LL( \x_{t+1},\y_{t+1}) + \LL(\x_{t+1},\y_{t+1}) - \LL( \x_{t+1},\y_{t}) - \eta_{t} \innerProd{M\x_{t+1}}{M\hat \x_{t+1}} \\
& \; \stackrel{(\star)}{=}\; (\LL(\x_{t+2},\y_{t+1}) - \LL(\x_{t+1},\y_{t+1})) + \eta_{t}\|M\x_{t+1}\|^2 - \eta_{t} \innerProd{M\x_{t+1}}{M\hat \x_{t+1}} \,,
\end{align*}
where the last inequality $(\star)$ is by definition of $\LL$ and because $\y_{t+1} - \y_{t} = \eta_tM\x_{t+1}$.
\endproof
We can now combine Lemma~\ref{lemma:dual_descent} and Lemma~\ref{lemma:primal_descent} with our technical result Cor.~\ref{cor:lower_bound_gradient} on the dual suboptimality to get our fundamental descent lemma only valid under Assumption~\ref{assump:main}.
\begin{lemma}[Fundamental descent Lemma]\label{lemma:fund_descent}
  Under Assumption~\ref{assump:main} we have that for all $t\geq 0$,
  \begin{equation*}
    \Delta_{t+1}-\Delta_{t}
   \leq \frac{2 \eta_{t}}{\lambda} \left( \LL(\x_{t+1},\y_{t+1})-\LL(\hat\x_{t+1},\y_{t+1}) \right) +  \LL(\x_{t+2},\y_{t+1})-\LL(\x_{t+1},\y_{t+1}) -\eta_{t} \frac{\alpha^2}{2L_\lambda D^2}\min\{ \Delta_{t+1}^{(d)}, \tfrac{L_\lambda D^2}{2}   \} \;,
  \end{equation*}
\end{lemma}
\proof
Combining Lemma~\ref{lemma:dual_descent} and Lemma~\ref{lemma:primal_descent} gives us,
\begin{align}
   \Delta_{t+1}-\Delta_{t}
    & =   [\Delta_{t+1}^{(p)}-\Delta_{t}^{(p)}]+[\Delta_{t+1}^{(d)}-\Delta_{t}^{(d)}]\notag\\
    &\leq  \eta_{t}\|M\x_{t+1}\|^2+ \LL(\x_{t+2},\y_{t+1})-\LL(\x_{t+1},\y_{t+1}) -2  \eta_{t} \innerProd{M\x_{t+1}}{M\hat \x_{t+1}} \notag \\
    & = \eta_{t}\|M\x_{t+1}-M\hat \x_{t+1}\|^2-\eta_{t}\|M \hat \x_{t+1}\|^2+ \LL(\x_{t+2},\y_{t+1})-\LL(\x_{t+1},\y_{t+1}) \label{eq:Mhat_x_t} \,.
\end{align}
Finally, from the ``strong convexity'' of $\LL(\cdot,\y_{t+1})$ respect to $M\x$ Prop.~\ref{prop:strong_conv_Mx} we obtain,
\begin{equation} \label{eq:fund_descent}
   \Delta_{t+1}-\Delta_{t}
   \leq \frac{2 \eta_{t}}{\lambda} \left( \LL(\x_{t+1},\y_{t+1})-\LL(\hat\x_{t+1},\y_{t+1}) \right) +  \LL(\x_{t+2},\y_{t+1})-\LL(\x_{t+1},\y_{t+1}) -\eta_{t}\|M\hat \x_{t+1}\|^2\;,
\end{equation}
where $\Delta_{t+1} := \Delta_{t+1}^{(p)}+\Delta_{t+1}^{(d)}$. Then we can use our fundamental technical result (Corollary~\eqref{cor:lower_bound_gradient}) relating the dual suboptimality and the norm of its gradient,
\begin{equation}
\|M\hat \x_{t+1}\|^2 \stackrel{Prop.\ref{prop:dual_smooth}}{=} \|\nabla d(\y_{t+1})\|^2
  \geq \frac{\alpha^2}{2L_\lambda D^2}\min\{ \Delta_{t+1}^{(d)}, \tfrac{L_\lambda D^2}{2}   \} \;,
 \end{equation}
 to get the desired lemma.
\endproof

The two following sections respectively deal with the proof of Theorem~\ref{thm:main_general} and Theorem~\ref{thm:main_polytope} they both start from our fundamental descent lemma (Lemma~\ref{lemma:fund_descent}).

\subsection{Proof of Theorem~\ref{thm:main_general}} %
\label{sub:proof_of_1_in_theorem_thm:main}
Let us first recall the setting and propose a detailed version of the first part of Thm.~\ref{thm:main_general}. The second part of Thm.~\ref{thm:main_general} is proposed in Corollary~\ref{cor:general_polytope_feasibility}.
\begin{reptheorem}{thm:main_general}
  If $\X$ is a compact convex set and $f$ is $L$-smooth,
using any algorithm with \emph{sublinear decrease}~\eqref{eq:def_subl_decrease} as inner loop in \FWAL\ \eqref{eq:fwal} and $\eta_t := \min\Big\{\frac{2}{\lambda}, \frac{\alpha^2}{2 \delta}\Big\}\frac{2}{t+2}$ then there exists a bounded $t_0 \geq 0$ such that,
  \begin{equation}
    \Delta_{t} \leq \min\left\{\frac{4\delta(t_0+2)}{t+2} ,\delta \right\}
   \quad \forall t \geq t_0
   \quad \text{and} \quad
   t_0 \leq \left( \tfrac{C}{\delta}+2 \right) \exp \left( \frac{\Delta_0 - \delta +2C}{2 \delta}  \right) \;.
\end{equation}
where $C := 8 \delta\max\Big(\frac{1}{4}, \frac{4 \delta^2}{\lambda^2 \alpha^4}\Big)\, .$ and $\delta := L_\lambda D^2$.

If we set $\eta_t = \min\Big\{\frac{2}{\lambda}, \frac{\alpha^2}{2\delta}\Big\} \frac{C}{\delta}$ for at least $t_0$ iterations and then $\eta_t := \min\Big\{\frac{2}{\lambda}, \frac{\alpha^2}{2 \delta}\Big\}\frac{2}{t+2}$ we get
\begin{equation}
    \Delta_{t} \leq \min\left\{\frac{4\delta(t_0+2)}{t+2} ,\delta \right\}
   \quad \forall t \geq t_0
   \quad \text{where} \quad
   t_0 = \max\Big\{1 + \frac{2(\Delta_0-\delta)C}{\delta^2}, \frac{C}{\delta}\Big\}
    \;.
\end{equation}
\end{reptheorem}

\proof
This proof will start from Lemma~\ref{lemma:fund_descent} and use the fact that if $\X$ is a general convex compact set, a usual Frank-Wolfe step with line search (Alg.~\ref{alg:FW}) produces a sublinear decrease~\eqref{eq:def_subl_decrease}. It leads to the following equation holding for any $\gamma \in [0,1]$,
\begin{equation}
   \Delta_{t+1}-\Delta_{t}
  \leq \left( \frac{2 \eta_{t}}{\lambda} - \gamma \right) \left( \LL(\x_{t+1},\y_{t+1})-\LL(\hat\x_{t+1},\y_{t+1}) \right) + \gamma^2 \frac{L_\lambda D^2}{2} - \eta_{t} \frac{\alpha^2}{2L_\lambda D^2}\min\{ \Delta_{t+1}^{(d)}, \tfrac{L_\lambda D^2}{2}  \} \;.
\end{equation}
Then for $\gamma = \frac{4 \eta_t}{\lambda}$ we get,
\begin{equation}
   \Delta_{t+1}-\Delta_{t}
  \leq  -\frac{2 \eta_{t}}{\lambda} \left( \LL(\x_{t+1},\y_{t+1})-\LL(\hat\x_{t+1},\y_{t+1}) \right) + \left( \frac{4 \eta_{t}}{\lambda} \right) ^2 \frac{L_\lambda D^2}{2} - \eta_{t} \frac{\alpha^2}{2L_\lambda D^2}\min\{ \Delta_{t+1}^{(d)}, \tfrac{L_\lambda D^2}{2} \} \;.
\end{equation}
Since we are doing line-search, we know that $ \LL(\x_{t+1},\y_{t+1}) \geq  \LL(\x_{t+2},\y_{t+1})$ implying that
\begin{equation}\label{eq:needed_following_cor}
   \Delta_{t+1}-\Delta_{t}
  \leq  -\frac{2 \eta_{t}}{\lambda} \Delta_{t+1}^{(p)} + \left( \frac{4 \eta_{t}}{\lambda} \right) ^2 \frac{L_\lambda D^2}{2} - \eta_{t} \frac{\alpha^2}{2L_\lambda D^2}\min\{ \Delta_{t+1}^{(d)}, \tfrac{L_\lambda D^2}{2} \} \;,
\end{equation}
In order to make appear $\Delta_{t+1}$ in the RHS, we will introduce
\begin{equation}
a = \min\Big\{\frac{2}{\lambda}, \frac{\alpha^2}{2L_\lambda D^2}\Big\} \,,
\end{equation}
this constant depends on $\lambda$ which is a hyperparameter.
It seems that $\lambda$ helps to scale the decrease of the primal with to the one of the dual.
\begin{equation}\label{eq:proof_sublin_pb_alternate_regime}
  \Delta_{t+1}-\Delta_{t}
  \leq  -a\eta_{t}\min\{ \Delta_{t+1}, \tfrac{L_\lambda D^2}{2} \} + \left( \frac{4 \eta_{t}}{\lambda} \right) ^2 \frac{L_\lambda D^2}{2}  \;.
\end{equation}

Then we have either that,
\begin{equation}
    \Delta_{t+1} - \Delta_{t} \leq  -a L_\lambda D^2 \eta_t/2 + \left( \frac{4 \eta_{t}}{\lambda} \right) ^2 \frac{L_\lambda D^2}{2} \;,
\end{equation}
giving a uniform (in time) decrease with a small enough constant step size $\eta_t$ or we have,
\begin{equation}
    \Delta_{t+1} - \Delta_{t} \leq - a\eta_t \Delta_{t+1} +\left( \frac{4\eta_{t}}{\lambda} \right) ^2 \frac{L_\lambda D^2}{2} \;,
\end{equation}
giving a usual Frank-Wolfe recurrence scheme leading to a sublinear decrease with a decreasing step size $\eta_t \sim 1/t$. It seems hard to get an adaptive step size since we cannot efficiently compute $\Delta_t$.
In order to tackle this problem we will consider an upper bound looser than \eqref{eq:proof_sublin_pb_alternate_regime} leading to a separation of the two regimes.
Let us introduce 
\begin{equation}
\bar \eta_t := a \eta_t\,, \; \delta := \tfrac{L_\lambda D^2}{2} 
\quad \text{and} \quad
C := 8 \delta\max\Big(\frac{1}{4}, \frac{4 \delta^2}{\lambda^2 \alpha^4}\Big)\, .
\end{equation}
Replacing $\eta_t$ with $\bar \eta_t$, we have that~\eqref{eq:proof_sublin_pb_alternate_regime} implies
\begin{equation}\label{eq:summarized_both_regimes}
  \Delta_{t+1}-\Delta_{t}
  \leq  -\bar\eta_{t}\min\{ \Delta_{t+1}, \delta \} +  \bar\eta_{t}^2 \frac{C}{2}  \;.
\end{equation}
\begin{lemma}
\label{lemma:sublin_decrese_subopt_frank_wolfe_scheme}
If there exists $t_0 > \frac{C}{\delta} -2$ such that $\Delta_{t_0} \leq \delta$ and if we set $\bar\eta_{t} = \frac{2}{2+t}$ then,
\begin{equation}\label{eq:hyp_rec_sublin}
   \Delta_{t} \leq \min\left\{\frac{4\delta(t_0+2)}{t+2} ,\delta \right\}
   \quad \forall t \geq t_0\;.
\end{equation}
\end{lemma}
\proof
For $t=t_0$ the result comes from the fact that we assumed that $\Delta_{t_0} \leq \delta$. By induction, let us assume that for a $t \geq  t_0, \;\Delta_{t} \leq \min\left\{\frac{4\delta(t_0+2)}{t+2} ,\delta \right\}$ then if $\Delta_{t+1}$ was greater than $\delta$, we would have obtained,
\begin{equation}
  \delta \leq \Delta_{t+1}
   \leq \Delta_t - \frac{2}{2+t}\delta + \left(\frac{2}{2+t} \right) ^2 \frac{C}{2} \leq \delta  - \frac{2}{2+t}\delta + \left(\frac{2}{2+t} \right) ^2 \frac{C}{2}\;,
\end{equation}
implying that,
\begin{equation}
  \delta \leq \frac{C}{2+t}  \qquad \text{and then} \qquad t \leq \frac{C}{\delta} - 2
\end{equation}
which contradicts the assumption $t > \frac{C}{\delta}-2$. Leading to $\Delta_t \leq \delta \,, \;\forall t \geq t_0$. 

Moreover, we have for all $t \geq t_0$,
\begin{align}
   \Delta_{t+1}
  & \leq \Delta_t - \frac{2}{2+t}\Delta_{t+1} + \left(\frac{2}{2+t} \right) ^2 \frac{C}{2} \\
  \frac{t+4}{t+2} \Delta_{t+1}
  & \leq \Delta_t +\left(\frac{2}{2+t} \right) ^2 \frac{C}{2} \\
  \Delta_{t+1}
  & \leq \frac{t+2}{t+4}\Delta_t + \frac{ 2 C}{(2+t)(t+4)} \\
  & \stackrel{(\star)}{\leq} \frac{t+2}{t+4}\frac{4 \delta (t_0+2)}{t+2} + \frac{ 2 C}{(t+2)(t+4)} \\
  & \leq \frac{4 \delta (t_0+2)}{t+3} \left[ \frac{t+3}{t+4} \left( 1 + \frac{ 1}{2(t+2)}\right)  \right] \;,
\end{align}
where ($\star$) is due to the induction hypothesis and the last inequality is due to the fact that $ \delta(t_0 + 2)\geq C$. Then, we just need to show that
\begin{equation}
  \left[ \frac{t+3}{t+4} \left( 1 + \frac{ 1}{2(t+2)}\right)  \right] \leq 1 \, ,\quad \forall t \geq 1\,.
\end{equation}
That is true because
\begin{align}
  & \left[ \frac{t+3}{t+4} \left( 1 + \frac{ 1}{2(t+2)}\right)  \right] \leq 1 \\
  \Leftrightarrow
  & \; (t+3)(t+\tfrac{5}{2}) \leq (2+t)(t+4) \\
  \Leftrightarrow
  & \;- \frac{1}{2}t + \frac{15}{2} \leq 8 \\
  \Leftrightarrow
  & \; \;t \geq 1\;.
\end{align}
\endproof
Now we have to show that in a finite number of iterations $t_0$ we can reach a point such that $\Delta_{t_0} \leq \delta$.

Let us assume that $\Delta_0 \geq \delta$, then we cannot initialize the recurrence~\eqref{eq:hyp_rec_sublin}. Instead we will show the following:
\begin{lemma}\label{lemma:upper_bound_t_0}
Let $(\Delta_t)_{t\in \N}$ a sequence such that $\Delta_{t+1}-\Delta_{t}
  \leq  -\bar\eta_{t}\min\{ \Delta_{t+1}, \delta \} +  \bar\eta_{t}^2 \frac{C}{2}\,,\; \forall t \in \N.\;$ We have that,
\begin{itemize}
   \item If $\bar\eta_t = \frac{\delta}{C}$, then there exists $t_0 \in \N$ such that,
\begin{equation}
     \Delta_{t_0} \leq \delta
     \,,\quad 
     \Delta_t \leq \delta \quad: \quad \forall\; t\geq t_0\;,
     \quad \text{and} \quad
     t_0 \leq 1+\frac{2(\Delta_0-\delta)C}{\delta^2}
     \;.
  \end{equation}
  \item If $\bar\eta_t = \frac{2}{2+t}$, then there exists $t_0 \geq \frac{C}{\delta} - 2$ such that,
\begin{equation}
     \Delta_{t_0} \leq \delta
     \quad \text{and} \quad
     t_0 \leq \left( \tfrac{C}{\delta} \right) \exp \left( \frac{\Delta_{\lfloor \frac{C}{\delta} -1\rfloor}  - \delta +2C}{2 \delta}  \right)
     \;.
  \end{equation}
 \end{itemize}
\end{lemma}
\proof
By contradiction, let us assume that $\Delta_t \geq \delta\,, \; \forall t$. Then \eqref{eq:summarized_both_regimes} gives,
\begin{equation}
  \Delta_{t+1}-\Delta_{t}
  \leq  -\bar\eta_{t} \delta  +  \bar\eta_{t}^2 \frac{C}{2}  \;.
\end{equation}
Then we would have for $\bar \eta_t = \frac{\delta}{C}$ that 
\begin{equation}
   \Delta_{t+1} \leq \Delta_{t} -  \frac{\delta^2}{2C}  
 \end{equation} 
Consequently we would have $\Delta_t <0$ at some point contradicting the fact that $\Delta_t$ is non negative.

For $\bar \eta_t = \frac{2}{2+t}$ we would have that,
\begin{equation}
 \infty = \delta \sum_{t=0}^\infty \bar \eta_t \leq \Delta_0 + \frac{C}{2}\sum_{t=0}^\infty \bar
 \eta_t^2 < \infty \;.
\end{equation}
giving a contradiction.

Thus, let us consider the smallest time $t_0$ such that $\Delta_{t_0} \leq \delta$.
\begin{itemize}
  \item
  If we set $\bar\eta_t = \frac{\delta}{C}$, we get for all $t < t_0 $
  \begin{equation}
    \Delta_{t+1}-\Delta_{t}
    \leq  -\frac{\delta^2}{2C}
  \end{equation}
  and then summing for $0 \leq t \leq t_0-2$
  \begin{equation}
    \delta - \Delta_0 \leq \Delta_{t_0-1} - \Delta_0 \leq - \frac{(t_0-1) \delta^2}{2C}\;,
  \end{equation}
  implying that
  \begin{equation}
    t_0 \leq 1 + \frac{2(\Delta_0-\delta)C}{\delta^2}
  \end{equation}
  then, let us show by recurrence that $\forall t\geq t_0, \; \Delta_t \leq \delta$. The result for $t=t_0$ is true by definition of $t_0$. Let us assume that it is true for a $t \geq t_0$, then if $\Delta_{t+1} \geq \delta$, \eqref{eq:summarized_both_regimes} gives us
  \begin{equation}
    \Delta_{t+1} \leq \Delta_t - \frac{\delta^2}{C} + \frac{\delta^2}{2C} \leq \Delta_t-  \frac{\delta^2}{2C}  < \delta \;,
  \end{equation}
  leading to a contradiction. Thus, we have that $\Delta_{t+1} \leq \delta$.
  \item If $\bar \eta_t = \frac{2}{2+t}$,  we want a $t_0 \geq \frac{C}{\delta} -2$ so if $\Delta_{\lfloor \frac{C}{\delta} -1\rfloor} \leq \delta$ we are done, otherwise
  \begin{equation}
    \delta \sum_{t=0}^{t_0 - 2} \frac{2}{2+t} \leq \Delta_0 - \delta + \frac{C}{2} \sum_{t=0}^\infty \frac{4}{(2+t)^2} \leq  \Delta_0 - \delta + 2C \Big(\frac{\pi^2}{6}-1\Big) \leq \Delta_{\lfloor \frac{C}{\delta} -1\rfloor}  - \delta + 2C  \;.
  \end{equation}
  Since $\sum_{t=\lfloor \frac{C}{\delta} -1 \rfloor}^{t_0} \frac{1}{t} \geq \ln(t_0) - \ln(\frac{C}{\delta})$ we get that
  \begin{equation}
    t_0 \leq \left( \tfrac{C}{\delta} \right) \exp \left( \frac{\Delta_{\lfloor \frac{C}{\delta} -1\rfloor}  - \delta +2C}{2 \delta}  \right) \;.
  \end{equation}
\end{itemize}
\endproof
Combining Lemma~\ref{lemma:sublin_decrese_subopt_frank_wolfe_scheme} and Lemma~\ref{lemma:upper_bound_t_0} with~\eqref{eq:summarized_both_regimes} we finally get Theorem~\ref{thm:main_general}.
\endproof

To sum up, we can either set $\bar \eta_t = \frac{\delta}{C}$ for a fixed number of iterations or we can use a decreasing step size leading to a very bad upper bound on $t_0$. 
Nevertheless this bound for the decreasing step size is very conservative and even if the best theoretical rates are given by a constant step size $\bar \eta_t$ for a number of iterations proportional to $\frac{C}{\delta}$ and then a sublinear step size $\bar \eta_t = \frac{2}{2+t}$, in practice, we can directly start with a decreasing step size.  

\begin{corollary}\label{cor:general_polytope_feasibility} Under the same assumption as Thm.~\ref{thm:main_general}. Let the $t_0 \in \N$ stated in Thm.~\ref{thm:main_general}, then for all $T\geq t_1\geq t_0$,
\begin{equation}
  \min_{t_1 \leq t \leq T}  \|M\hat \x_{t+1}\|^2
   \leq \frac{2\Theta}{T-t_1+1}
   \quad \text{and} \quad
   \min_{t_1 \leq t \leq T}  \|M\x_{t+1}\|^2
   \leq \frac{8\Theta}{T-t_1+1} \;.
\end{equation}
where $\Theta :=  \frac{8 \delta(t_0+2)}{a} + \frac{16a \delta}{\lambda^2}  $. Moreover, if $f$ is $\mu$-strongly convex we have for all $t_1 \geq 8t_0+15$,
   \begin{equation}
     \min_{t_1 \leq t \leq T} \|\x - \x^*\|^2 \leq \frac{4K}{\mu} \left[\frac{\lambda}{2} + \frac{2^8\delta^2}{\alpha^4 \lambda} \right]\frac{\Theta}{T-t_1+1}\;.
   \end{equation}
\end{corollary}
\proof This proof follows the same idea as the proof of~\citep[Thm C.3]{lacostejulien:hal-00720158}. Since we are working with different quantities and that the rates are slightly different from the ones provided in~\citep{lacostejulien:hal-00720158} we will provide a complete proof of this result.
We start from the fundamental descent lemma~\eqref{eq:fund_descent}. We use the fact that a usual Frank-Wolfe step produces a sublinear decrease~\eqref{par:sublinear_decrease_} that we specify for $\gamma = \frac{4\eta_t}{\lambda}$ to get a similar equation as~\eqref{eq:needed_following_cor},
\begin{equation}
   \Delta_{t+1}-\Delta_{t}
  \leq  -\frac{2 \eta_{t}}{\lambda} \left( \LL(\x_{t+1},\y_{t+1})-\LL(\hat\x_{t+1},\y_{t+1}) \right) + \left( \frac{4 \eta_{t}}{\lambda} \right) ^2 \frac{\delta}{2} - \eta_{t} \|M\hat \x_{t+1}\|^2 \;,
\end{equation}
noting $\delta := L_\lambda D^2$. Then introducing $h_{t+1} := \left( \LL(\x_{t+1},\y_{t+1})-\LL(\hat\x_{t+1},\y_{t+1}) \right)$, (note that because of the line search $\Delta^{(p)}_{t+1} \geq h_{t+1} \geq 0$) it leads to
\begin{equation}\label{eq:feasability_upper_bounded}
   \left( \frac{2}{\lambda}h_{t+1} +  \|M\hat \x_{t+1}\|^2  \right)
  \leq  \frac{\Delta_{t}-\Delta_{t+1}}{\eta_t}  +\eta_{t} \left( \frac{4 }{\lambda} \right) ^2 \frac{\delta}{2} \;.
\end{equation}
which is similar equation as~\citep[Eq.(22)]{lacostejulien:hal-00720158}. 
Let $t_1 \geq t_0$ and $\{w_t\}_{t_1}^T$ be a sequence of positive weights.  Let $\rho_t := \nicefrac{w_t}{\sum_{t=t_1}^T w_t}$ be the associated normalized weights. The convex combination of~\eqref{eq:feasability_upper_bounded} give us,
\begin{align}
  \sum_{t= t_1}^T \rho_t \left( \frac{2}{\lambda}h_{t+1} +  \|M\hat \x_{t+1}\|^2  \right)
  & \leq  \sum_{t= t_1}^T \rho_t \frac{\Delta_{t}-\Delta_{t+1}}{\eta_t}  + \sum_{t= t_1}^T \rho_t\eta_{t} \left( \frac{4}{\lambda} \right) ^2 \frac{\delta}{2} \notag\\
  & = \frac{\rho_{t_1}}{\eta_{t_1}} \Delta_{t_1} - \frac{\rho_{T}}{\eta_{T}} \Delta_{T+1} + \sum_{t= t_1}^{T-1} \Delta_{t} \left(\frac{\rho_{t+1}}{\eta_{t+1}} - \frac{\rho_t}{\eta_t}  \right) + \sum_{t= t_1}^T \rho_t\eta_{t} \left( \frac{4 }{\lambda} \right) ^2 \frac{\delta}{2} \notag\\
  & \leq \frac{\rho_{t_1}}{\eta_{t_1}} \Delta_{t_1} +  \sum_{t= t_1}^{T-1} \Delta_{t} \left(\frac{\rho_{t+1}}{\eta_{t+1}} - \frac{\rho_t}{\eta_t}  \right) + \sum_{t= t_1}^T \rho_t\eta_{t} \left( \frac{4 }{\lambda} \right) ^2 \frac{\delta}{2} \;.
\end{align}
We can now use a \emph{weighted average} such as $w_t = t- t_1$. This kind of average leads to
\begin{equation}
  \frac{\rho_{t+1}}{\eta_{t+1}} - \frac{\rho_t}{\eta_t} = \frac{ (t-t_1+1)(t+3) - (t-t_1)(t+2) }{a(T-t_1)(T-t_1+1)} = \frac{2t - t_1 + 3}{a(T-t_1)(T-t_1+1)}\;.
\end{equation}
where $\eta_t := \min\Big\{\frac{2}{\lambda}, \frac{\alpha^2}{4 \delta}\Big\}\frac{2}{t+2} =  a\frac{2}{t+2}$.
Then we can plug that $\Delta_t \leq \min\left\{\frac{4 \delta (t_0+2)}{t+2},\delta\right\}, \; \forall t \geq t_1 \geq t_0$ to get,
\begin{align}
  \sum_{t= t_1}^T \rho_t \left( \frac{2}{\lambda}h_{t+1} +  \|M\hat \x_{t+1}\|^2  \right)
  & \leq \frac{\delta \rho_{t_1}}{\eta_{t_1}} +  \sum_{t= t_1}^{T-1} \frac{4 \delta (t_0+2)}{t+2} \left(\frac{\rho_{t+1}}{\eta_{t+1}} - \frac{\rho_t}{\eta_t}  \right) + \sum_{t= t_1}^T \rho_t\eta_{t} \left( \frac{4 }{\lambda} \right) ^2 \frac{\delta}{2} \\
  & \leq \frac{2}{(T-t_1)(T-t_1+1)}\left[ \sum_{t= t_1}^{T-1} \tfrac{4 \delta (t_0+2)(2t-t_1+3)}{a(t+2)}  + \sum_{t= t_1}^T \tfrac{2a(t-t_1)}{2+t} \left( \tfrac{4 }{\lambda} \right) ^2 \tfrac{\delta}{2}\right] \\
  & \leq \frac{2}{T-t_1+1} \left[ \frac{8 \delta(t_0+2)}{a} + \frac{16a \delta}{\lambda^2} \right]
\end{align}
Then,
\begin{equation}
  \min_{t_1 \leq t \leq T}\left[ \frac{2}{\lambda}h_{t+1} +  \|M\hat \x_{t+1}\|^2 \right]
   \leq \frac{2}{T-t_1+1} \left[ \frac{8 \delta(t_0+2)}{a} + \frac{16a \delta}{\lambda^2} \right]\;.
\end{equation}
To upper bound $\|M\x_{t+1}\|^2$ the idea is to combine the previous equation with $\|M\x_{t+1} - M \hat\x_{t}\|^2 \leq \frac{2}{\lambda} \left(\LL(\x_{t+1},\y_{t}) - \LL(\hat \x_{t},\y_{t})\right) =:\frac{2}{\lambda} \Delta^{(p)}_{t} \leq \frac{2}{\lambda} h_t$ (Prop.~\ref{prop:strong_conv_Mx} plus the fact that we perform line search) giving,
\begin{equation}\label{eq:ineg_feas_sub_dual_gradient}
  \|M \x_{t+1} \|^2 \leq \frac{8}{\lambda} h_t + 4 \|M \hat \x_{t}\|^2
\end{equation}
  \begin{equation}
     \min_{t_1 +1 \leq t \leq T+1}\|M \x_{t+1}\|^2 \leq 4 \min_{t_1+1 \leq t \leq T+1} \left(\| M \hat \x_{t}\|^2 + \frac{2}{\lambda}h_{t} \right) \leq \frac{8}{T-t_1+1} \left[ \frac{8 \delta(t_0+2)}{a} + \frac{16a \delta}{\lambda^2} \right]
  \end{equation}
  If $f$ is $\mu$-strongly convex we can use Prop.~\ref{prop:delta_bounds_input} to get,
  \begin{equation}
   \frac{\mu}{2}\|\x_{t+1}-\x^*\|^2 \leq \Delta_t^{(p)} + \frac{2\delta}{\alpha^2} \|M\x_{t+1}\| \|M\hat \x_{t+1}\|\,, \; \forall t \in \N \;;\; \dist(\y_t,\Y^*) \leq \frac{\delta}{\alpha}\;.
\end{equation}
In order to show that at some point we have $\dist(\y_t,\Y^*) \leq \frac{\delta}{\alpha}$ we will use Thm.~\ref{thm:dual_function_PL} and \eqref{eq:hyp_rec_sublin} to get,
\begin{equation}
  \frac{4 \delta (t_0+2)}{t+2}
  \geq \Delta^{p}_{t}
  \geq \frac{1}{2\delta}\min \left\{
        \alpha^2\dist(\y_t,\Y^*)^2, \alpha \delta\dist(\y_t,\Y^*)
  \right\} \;, \; \forall t \geq t_0 \,,
\end{equation}
Then for all $t \geq t_0$ such that $\dist(\y_t,\Y^*) > \frac{\delta}{\alpha}$ we have that,
\begin{equation}
  \frac{8 \delta (t_0+2)}{t+2} \geq  \alpha \dist(\y_t,\Y^*)
\end{equation}
implying that for $t\geq 8(t_0+2) -2 = 8t_0 + 14$ we have that $\alpha\dist(\y_t,\Y^*) \leq \delta$ and then,
\begin{align}
   \frac{\mu}{2}\|\x_{t+1}-\x^*\|^2
   &\leq \Delta_t^{(p)} + \frac{2\delta}{\alpha^2} \|M\x_{t+1}\| \|M\hat \x_{t}\| \\
   & \stackrel{\eqref{eq:ineg_feas_sub_dual_gradient}}{\leq} h_{t} + \frac{2\delta}{\alpha^2}  4 \sqrt{2}\|M\hat \x_{t}\| \sqrt{\frac{2}{\lambda} h_{t} + \|M \hat \x_{t}\|^2 } \\
   & \leq \frac{\lambda}{2} \|M\hat \x_t\|^2 + h_t + \frac{2^8\delta^2}{\alpha^4 \lambda} \left( \frac{2}{\lambda} h_{t} + \|M \hat \x_{t}\|^2 \right) \\
   & \leq \left( \frac{\lambda}{2} + \frac{2^8\delta^2}{\alpha^4 \lambda}  \right) \left( \frac{2}{\lambda} h_{t} + \|M \hat \x_{t}\|^2 \right)\;.
\end{align}
It then implies that for $t_1 \geq 8t_0+14$,
\begin{align}
     \min_{t_1+ 1 \leq t \leq T+1} \|\x_{t+1} - \x^*\|^2
     &\leq \frac{2}{\mu}  \left( \frac{\lambda}{2} + \frac{2^8\delta^2}{\alpha^4 \lambda}  \right) \min_{t_1+ 1\leq t \leq T+1}\left( \frac{2}{\lambda} h_{t} + \|M \hat \x_{t}\|^2 \right) \\
     & \leq \frac{2}{\mu}  \left( \frac{\lambda}{2} + \frac{2^8\delta^2}{\alpha^4 \lambda}  \right)\frac{2}{T-t_1+ 1} \left[ \frac{8 \delta(t_0+2)}{a} + \frac{16a \delta}{\lambda^2} \right]
     \;.
   \end{align}
\endproof
\subsection{Proof of Theorem~\ref{thm:main_polytope}} %
\label{sub:proof_of_theorem_thm:main_polytope}
This proof starts with the fundamental descent lemma (Lemma~\ref{lemma:fund_descent}). It uses the fact that if $\X$ is a polytope and if we use an algorithm with a geometric decrease~\eqref{eq:def_geom_decrease} such as Alg.~\ref{alg:AFW} then with a small enough constant step size $\eta_{t}$ we can upper bound the decrease of $\Delta_{t+1} - \Delta_{t}$.

\begin{lemma} \label{lemma:geom_decrease_delta}
Under assumptions of theorem~\ref{thm:main_polytope}, we have 
\begin{equation}
  \Delta_{t+1}-\Delta_{t}
  \leq -\frac{\rho_A}{2 } \Delta_{t+1}^{(p)} -\frac{\lambda \alpha^2 \rho_A}{8 L_\lambda D^2} \min\{\Delta_{t+1}^{(d)},\tfrac{L_\lambda D^2}{2} \} 
  \quad \forall t \geq 1\;.
\end{equation}
\end{lemma}
\proof
To prove Lemma~\ref{lemma:geom_decrease_delta},
we start from Lemma~\ref{lemma:fund_descent} to obtain
\begin{align}
\Delta_{t+1}-\Delta_{t}
  &\leq   \LL(\x_{t+2},\y_{t+1})-\LL(\x_{t+1},\y_{t+1}) + \frac{2 \eta_{t}}{\lambda} \left( \LL(\x_{t+1},\y_{t+1})-\LL(\hat\x_{t+1},\y_{t+1}) \right) -\eta_{t} \frac{\alpha^2}{2L_\lambda D^2} \min\{\Delta_{t+1}^{(d)},\tfrac{L_\lambda D^2}{2} \} \;, \nonumber \\
  &\!\!\stackrel{\eqref{eq:def_geom_decrease}}{\leq}   \left( \frac{2 \eta_{t}}{\lambda} - \rho_A \right)  \left( \LL(\x_{t+1},\y_{t+1}) - \LL(\hat \x_{t+1},\y_{t+1}) \right) -\eta_{t} \frac{\alpha^2}{2L_\lambda D^2} \min\{\Delta_{t+1}^{(d)},\tfrac{L_\lambda D^2}{2} \} \;.
\label{eq:estimate}
\end{align}
Now we can choose $\eta_{t} = \frac{\lambda \cdot \rho_A}{4}$ giving us Lemma~\ref{lemma:geom_decrease_delta}.
\endproof
From this lemma we can deduce a constant decrease for a finite number of step and eventually a geometric decrease.
\begin{lemma} \label{lemma:geom_conv_delta}
Under the assumptions of Theorem~\ref{thm:main_polytope}, for all $\lambda > 0 $, if we set $\eta_t = \frac{\lambda \rho_A}{4}$ for finite number of steps $t_0$, then the quantity $\Delta_t$ decreases by a uniform amount as,
\begin{equation}
    \Delta_{t+1}-\Delta_{t}
  \leq  -\frac{\lambda \alpha^2 \rho_A}{16}
  \quad \text{where} \quad
  t_0(\Delta_0) \leq 1 + \frac{16\Delta_0 - 8L_\lambda D^2}{\lambda \alpha^2}\;.
\end{equation}
Otherwise, $\Delta_t$ decrease geometrically as,
\begin{equation}
    \Delta_{t+1}
  \leq \frac{1}{1+\rho}\Delta_{t}
  \quad \text{where} \quad
  \rho := \frac{\rho_A}{2} \min \left\{ 1, \frac{\lambda \alpha^2}{4 L_\lambda D^2}\right\} \;.
\end{equation}
\end{lemma}
\proof
We start from Lemma~\ref{lemma:geom_decrease_delta}, if $2\Delta_{t+1}^{(d)} \geq L_\lambda D^2$,
\begin{equation}\label{eq:uniform_decrease_delta}
  \Delta_{t+1} -\Delta_t \leq  -\frac{\lambda \alpha^2 \rho_A}{16} \;,
\end{equation}
and otherwise,
\begin{equation}
  \Delta_{t+1}-\Delta_{t}
  \leq -\frac{\rho_A}{2 } \Delta_{t+1}^{(p)} -\frac{\lambda \alpha^2 \rho_A}{8 L_\lambda D^2}\Delta_{t+1}^{(d)}
  \leq -\rho \Delta_{t+1}\;.
\end{equation}
Our goal is then just to show the upper bound on $t_0$. First let us notice that $(\Delta_t)$ is decreasing then this non negative sequence cannot decrease by a uniform amount an infinite number of time, then we can sum for $t = 1,\ldots,t_0-1$ such that \eqref{eq:uniform_decrease_delta} holds to get,
\begin{equation}
  \tfrac{L_\lambda D^2}{2}  - \Delta_0 \leq\Delta^{(d)}_{t_0-1}-\Delta_0  \leq \sum_{t=0}^{t_0-2} \Delta_{t+1}-\Delta_{t} \leq -(t_0-1)  \frac{\lambda \alpha^2 \rho_A}{16}
\end{equation}
\endproof
One can deduce several convergence properties from this lemma which are compiled in Theorem~\ref{thm:main_polytope}.
\begin{corollary}[Extended Theorem~\ref{thm:main_polytope}]\label{cor:main_polytope}
Under the assumptions of Theorem~\ref{thm:main_polytope},
there exist $t_0 \leq1+\frac{16\Delta_0 - 8L_\lambda D^2}{\lambda \alpha^2}$ such that
for all $t\geq t_0$ we have the following properties,
\begin{enumerate}
  \item The gap decreases linearly,
  \begin{equation}
    \Delta_{t} \leq \frac{L_\lambda D^2}{2(1+\rho)^{t-t_0}} \;.
  \end{equation}
  \item The sequences of feasibility violations at points $\x_{t+1}$ and $\hat \x_{t+1}$ decrease linearly,
  \begin{equation}
    \|M\hat \x_{t+1}\|^2 \leq  \frac{2}{\lambda \cdot \rho_A} \frac{L_\lambda D^2}{(1+\rho)^{t-t_0}} \quad \text{and} \quad \|M\x_{t+1}\|^2 \leq  \frac{8}{\lambda \cdot \rho_A} \frac{L_\lambda D^2}{(1+\rho)^{t-t_0}}, \quad
  \end{equation}
\end{enumerate}
where $\rho : =  \frac{\rho_A}{2} \min \left\{ 1, \frac{\lambda \alpha^2}{8 L_\lambda D^2}\right\}$.
  \item Finally, if $f$ is $\mu_f$-strongly convex, the distance of the current point to the optimal set vanishes as,
  \begin{equation}
    \|\x_{t+1} -\x^*\|^2  \leq \frac{L_\lambda D^2(\sqrt{2}+1)}{2\mu_f(\sqrt{1+\rho})^{t-t_0}} + \frac{O(1)}{(1+\rho)^{t-t_0}}\;.
  \end{equation}
\end{corollary}

\proof

To prove the first statement let us start from Lemma~\ref{lemma:geom_decrease_delta}, for all $t \geq t_0$,
\begin{equation}
  [\Delta_{t+1}^{(p)}+\Delta_{t+1}^{(d)}]-[\Delta_{t}^{(p)}+\Delta_{t}^{(d)}]
  \leq -\frac{\rho_A}{2 } \left( \LL(\x_{t+1},\y_{t+1}) - \LL(\hat \x_{t+1},\y_{t+1}) \right) -\frac{\lambda \cdot \rho_A}{4}\|M\hat \x_{t+1}\|^2 \;,
\end{equation}
  leading us directly to
  \begin{equation}
    \frac{\lambda \cdot \rho_A}{4}\|M \hat \x_{t+1}\|^2 \leq \Delta_{t}
    \quad \text{and} \quad
    \frac{\rho_A}{2 } \left( \LL(\x_{t+1},\y_{t+1}) - \LL(\hat \x_{t+1},\y_{t+1}) \right) \leq \Delta_{t} \;.
  \end{equation}
  To upper bound $\|M\x_{t+1}\|^2$ the idea is to combine the two previous equations with $\|M\x_{t+1} - M \hat\x_{t+1}\|^2 \leq \frac{2}{\lambda} \left(\LL(\x_{t+1},\y_{t+1}) - \LL(\hat \x_{t+1},\y_{t+1})\right)$ (Prop.~\ref{prop:strong_conv_Mx}) giving,
  \begin{equation}
    \|M \x_{t+1}\|^2 \leq 2 \| M \hat \x_{t+1}\|^2 + 2\|M\x_{t+1} - M \hat\x_{t+1}\|^2 \leq \frac{16}{\lambda \rho_A} \Delta_{t} \;.
  \end{equation}

The last statement directly follows from the fact $(\Delta_{t+1})$ decreases linearly (Lemma~\ref{lemma:geom_conv_delta}) and the fact that one can upper bound $\|\x_{t+1} - \x^*\|^2$ with the primal and dual suboptimalities (Proposition~\ref{prop:delta_bounds_input}),
\begin{equation}
   \frac{\mu}{2}\|\x_{t+1}-\x^*\|^2 \leq \Delta_t^{(p)} - \Delta_t^{(d)} + \max \left( 2\Delta_t^{(d)}, \sqrt{2L_\lambda D^2 \Delta_t^{(d)}} \right)\;.
\end{equation}
Then, it easily follows that
\begin{equation}
  \frac{\mu}{2}\|\x_{t+1}-\x^*\|^2 \leq \frac{(\sqrt{2}+1)L_\lambda D^2}{(\sqrt{1+\rho})^{t-t_0}} + \frac{O(1)}{(1+\rho)^{t-t_0}} \;.
\end{equation}
\endproof
\end{document}